\documentclass{article}
\usepackage{amsmath,epsfig,latexsym,amssymb,cite,graphicx,a4}
\newlength{\indentedwidth} \newdimen\mathindent
\indentedwidth=\mathindent

\newtheorem{proposition}{Proposition}[section]

\newtheorem{definition}{Definition}[section]

\usepackage{epsf,graphics}
\DeclareMathAlphabet{\mathpzc}{OT1}{pzc}{m}{it}
\begin{document}
\vskip 0.5cm
\begin{center}
{\Large \bf Lax Operator for the Quantised Orthosymplectic
Superalgebra $U_q[osp(m|n)]$}
\end{center}
\vskip 0.8cm
\centerline{K.A. Dancer,%
\footnote{\tt dancer@maths.uq.edu.au} M.D. Gould%
\footnote{\tt mdg@maths.uq.edu.au} and J. Links%
\footnote{\tt jrl@maths.uq.edu.au}}
\vskip 0.9cm \centerline{\sl\small Centre for Mathematical Physics,
School of Physical Sciences, } \centerline{\sl\small The University of
Queensland, Brisbane 4072, Australia.}
\vskip 0.9cm

\begin{abstract}
\noindent Representations of quantum superalgebras provide a natural
framework in which to model supersymmetric quantum systems.  Each quantum
superalgebra, belonging to the class of quasi-triangular Hopf
superalgebras, contains a universal $R$-matrix which automatically
satisfies the Yang--Baxter equation.  Applying the vector
representation $\pi$, which acts on the vector module $V$, to the
left-hand side of a universal $R$-matrix gives a Lax operator.  In
this Communication a Lax operator is constructed for the quantised
orthosymplectic superalgebras $U_q[osp(m|n)]$ for all $m > 2, n \geq
0$ where $n$ is even.  This can then be used to find a solution to the
Yang--Baxter equation acting on $V\otimes V\otimes W$, where $W$ is an
arbitrary $U_q[osp(m|n)]$ module.  The case $W=V$ is studied as an
example.
\end{abstract}


\section{Introduction}

Solutions to the (spectral parameter dependent) Yang-Baxter equation
(YBE) lie at the core of the Quantum Inverse Scattering Method for the
construction of integrable quantum systems, and underlie the
applicability of Bethe ansatz methods for deriving their exact
solutions (e.g. see \cite{Baxter,Jimbo1}).  Integrable systems with
both bosonic and fermionic degrees of freedom may be constructed
through a $\mathbb{Z}_2$-graded, or \emph{supersymmetric}, analogue of
the YBE. Several examples exist for models constructed in this manner
which can be used to describe systems of strongly correlated
electronic \cite{fk1,ek,eks,fk2,bglz,rm,glzt,mr} and also integrable
supersymmetric field theories \cite{Saleur1,Saleur2,Saleur3}.

A systematic approach to solve the YBE is provided by the Quantum
Double construction \cite{Drinfeld} (see \cite{gzb} for the
$\mathbb{Z}_2$-graded extension). This gives a prescription for
embedding a Hopf algebra and its dual into a quasi-triangular Hopf
algebra $A$, which possesses a \emph{universal} $R$-matrix
$\mathcal{R}\in\, A \otimes A$ providing an algebraic solution for the
YBE. Each tensor product matrix representation of the quasi-triangular
Hopf algebra then yields a matrix solution for the YBE.  The most
well-known examples of quasi-triangular Hopf superalgebras are the
quantum superalgebras \cite{Yamane} (denoted $U_q[g]$ where $g$ is a
Lie superalgebra), which are deformations of the universal enveloping
algebras of Lie superalgebras. In cases where $g$ is affine, loop
representations of $U_q[g]$ provide supersymmetric matrix solutions of
the YBE dependent on a spectral parameter $u$. The limit $u\rightarrow
\infty$ gives a solution associated with a non-affine subsuperalgebra
of $g$. It is possible however to start with a solution for the
non-affine case and then introduce a spectral parameter, a process
known as \emph{Baxterisation}. A systematic way to perform
Baxterisation in the case of quantum superalgebras, through the use of
\emph{tensor product graph methods}, is described in
\cite{dglz,GouldZhang00}.  This result then reduces the problem to
determining the explicit form for $R$-matrices associated with
non-affine quantum superalgebras.

In practice however, it is still a difficult technical challenge to
explicitly compute $\mathcal{R}$ for a given non-affine $U_q[g]$.  Our
approach here is to simplify the problem by instead looking to
determine the Lax operator $R=(\pi \otimes {\rm id})\mathcal{R}$ where
$\pi$ denotes the vector representation of $U_q[g]$. Explicit forms
for the Lax operator are known in the case of $U_q[gl(m)]$
\cite{Jimbo2} and $U_q[gl(m|n)]$ \cite{Zhang2} but even in the
non-graded cases of $U_q[o(m)]$ and $U_q[sp(n)]$ the Lax operator is
not known.  In this Communication a Lax operator is constructed for
the quantised orthosymplectic superalgebras $U_q[osp(m|n)]$ for all $m
> 2, n \geq 0$ where $n$ is even. The case with $n=0$ corresponds to 
 constructing a Lax operator for $U_q[o(m)]$.  For technical reasons, 
the cases
$m=1$ and $m=2$ do not fit into the general framework that we will
employ here. This issue will be discussed in more detail in Sect.
\ref{q-def}. We mention here that in the case $m=1$ the Lax operator
can be deduced from an isomorphism given in \cite{Zhang}, while the
Lax operator for the case $m=2$ is derived elsewhere \cite{CLax}.  One
immediate application of the Lax operator is that it affords a means to
construct the Casimir invariants of the associated quantum
superalgebra. This construction for $U_q[osp(m|n)]$, as well as a
derivation of the eigenvalues of the Casimir invariants when acting on
irreducible highest weight modules, can be found in \cite{Casimir}.

We begin in Sect. 2 with a construction for $osp(m|n)$ which is
necessary to establish our notational conventions, and in particular
the choice of root system. Sect. 3 discusses the $q$-deformation of
the universal enveloping algebra and the associated quasi-triangular
Hopf superalgebra structure. In Sect. 4 we describe the construction
of the Lax operator by choosing a particular ansatz and imposing that
the Lax operator intertwines the co-product. It will be shown that,
for this choice of ansatz, there is a unique solution satsifying the
intertwining property. We will further show that this solution
satisfies all other properties which follow from the quasi-triangular
structure of the quantum orthosymplectic superalgebras.  Through use
of the Lax operator we derive the $R$-matrix for the vector
representation in Sect. 5, and concluding remarks are given in
Sect. 6.

\section{The Construction of $osp(m|n)$}

To construct the quantised orthosymplectic superalgebra
$U_q[osp(m|n)]$ we closely follow the method used in
\cite{GouldZhang00,GouldZhang99}, but provide details here in order to
establish our notational conventions.  We begin by developing
$osp(m|n)$ as a graded subalgebra of $gl(m|n)$.  The enveloping
algebra of $osp(m|n)$ is then deformed to yield $U_q[osp(m|n)]$, which
reduces to the original enveloping superalgebra as $q \rightarrow 1$.
The construction of $osp(m|n)$ starts with the standard generators
$e^a_b$ of $gl(m|n)$, the $(m+n) \times (m+n)$-dimensional general
linear superalgebra, whose even part is given by $gl(m) \oplus gl(n)$.
The commutator for a $\mathbb{Z}_2$-graded algebra satisfies the
relation

\begin{equation*}
[A,B] = - (-1)^{[A][B]}[B,A],
\end{equation*}
where $A,B$ are homogeneous operators and $[A] \in \mathbb{Z}_2$ is
the grading of $A$.  In particular, the generators of $gl(m|n)$
satisfy the graded commutation relations

\begin{equation*}
[e^{a}_{b},e^{c}_{d}] = \delta^{c}_{b} e^{a}_{d} -
(-1)^{([a]+[b])([c]+[d])} \delta^{a}_{d} e^{c}_{b}
\end{equation*}
where

\begin{equation*}
[a] =
   \begin{cases}
     0, \qquad \; a=i, & 1 \leq i \leq m, \\ 1, \qquad \; a = \mu, & 1
     \leq \mu \leq n.
   \end {cases}
\end{equation*}

Throughout this Communication we use Greek letters $\mu, \nu$ etc. to
denote odd indices and Latin letters $i, j$ etc. for even indices.  If
the grading is unknown, the usual $a, b, c$ etc. are used.  Which
convention applies will be clear from the context.  We will only ever
consider the homogeneous elements, but all results can be extended to
the inhomogeneous elements by linearity.

The orthosymplectic superalgebra $osp(m|n)$ is a subsuperalgebra of
$gl(m|n)$ with even part equal to $o(m) \oplus sp(n)$, where $o(m)$ is
the orthogonal Lie algebra of rank $\lfloor \frac{m}{2} \rfloor$ and $sp(n)$ 
is the symplectic
Lie algebra of rank $\frac{n}{2}$.  The latter only exists if $n$ is even, so
we set $n=2k$.  We also set $l = \lfloor \frac{m}{2} \rfloor$, so
$m=2l$ or $m=2l+1$.  To construct $osp(m|n)$ we require an even
non-degenerate supersymmetric metric $g_{ab}$.  Any can be used, but
for the sake of simplicity we choose $g_{ab} = \xi_{a}
\delta^a_{\overline{b}}$, with inverse metric $g^{ba} = \xi_{b}
\delta^a_{\overline{b}}$. Here

\begin{equation*}
\overline{a}=
  \begin{cases}
     m + 1 - a, & [a]=0, \\ n + 1 - a, & [a]=1,
  \end {cases}
\qquad \text{and} \quad \xi_{a} =
  \begin{cases}
     1, & [a] = 0, \\ (-1)^a, & [a]=1.
  \end{cases}
\end{equation*}

\

The $\mathbb{Z}_2$-graded subalgebra $osp(m|n)$ actually arises
naturally from considering the automorphism $\omega$ of $gl(m|n=2k)$
given by:

\begin{equation*}
\omega (e^a_b) = - (-1)^{[a]([a]+[b])} \xi_a \xi_b
e^{\overline{b}}_{\overline{a}}.
\end{equation*}

\noindent This is clearly of degree 2, with eigenvalues $\pm 1$, so it gives a
decomposition of $gl(m|n)$:

\begin{equation*}
gl(m|n) = \mathcal{S} \oplus \mathcal{T}, \text{ with }
[\mathcal{S},\mathcal{S}] \subset \mathcal{S},\;
[\mathcal{T},\mathcal{T}] \subset \mathcal{S} \text{ and }
[\mathcal{S},\mathcal{T}] \subset \mathcal{T},
\end{equation*}
where

\begin{alignat*}{2}
&\omega (x) = x && \forall x \in \mathcal{S}, \\ &\omega (x) = -x
&\qquad & \forall x \in \mathcal{T}.
\end{alignat*}
Here $\mathcal{T}$ is generated by operators

\begin{equation*}
T_{ab} = g_{ac} e^c_b + (-1)^{[a][b]} g_{bc} e^c_a = (-1)^{[a][b]}
T_{ba},
\end{equation*}
while $\mathcal{S}$ is generated by

\begin{equation*}
\sigma_{ab} = g_{ac} e^c_b - (-1)^{[a][b]} g_{bc} e^c_a =
-(-1)^{[a][b]} \sigma_{ba}.
\end{equation*}

\noindent The fixed-point $\mathbb{Z}_2$-graded subalgebra
$\mathcal{S}$ is the orthosymplectic superalgebra $osp(m|n)$, with the
operators $\sigma_{ab}$ providing a basis.  They satisfy the
commutation relations

\begin{multline*}
[\sigma_{ab}, \sigma_{cd}] = g_{cb} \sigma_{ad} -
  (-1)^{([a]+[b])([c]+[d])} g_{ad} \sigma_{cb} \\ - (-1)^{[c][d]}
  \bigl( g_{db} \sigma_{ac} - (-1)^{([a]+[b])([c]+[d])} g_{ac}
  \sigma_{db} \bigr).
\end{multline*}

As a more convenient basis for $osp(m|n)$ we choose the set of
Cartan-Weyl generators, given by:

\begin{align}
\sigma^a_b &= g^{ac} \sigma_{cb} \notag \\ &= e^a_b -
(-1)^{[a]([a]+[b])} \xi_a \xi_b
e^{\overline{b}}_{\overline{a}}. \label{CW}
\end{align}

\noindent Then the Cartan subalgebra $H$ is generated by the diagonal
operators

\begin{equation*}
\sigma^a_a = e^a_a - e^{\overline{a}}_{\overline{a}},
\end{equation*}
which satisfy $$[\sigma^a_a, \sigma^b_b]= 0, \qquad \forall a,b.$$

\

As a weight system, we take the set $\{ \varepsilon_i,\; 1 \leq i \leq
m \} \cup \{ \delta_\mu,\; 1 \leq \mu \leq n \}$, where
\mbox{$\varepsilon_{\overline{i}} = - \varepsilon_i$} and
\mbox{$\delta_{\overline{\mu}} = - \delta_\mu$}.  Conveniently, when
$m=2l+1$ this implies $\varepsilon_{l+1} = -\varepsilon_{l+1} = 0$.
Acting on these weights, we have the invariant bilinear form defined
by:

\begin{equation*}
(\varepsilon_i, \varepsilon_j) = \delta^i_j, \quad (\delta_\mu,
\delta_\nu) = -\delta^\mu_\nu, \quad (\varepsilon_i, \delta_\mu) = 0,
\qquad 1 \leq i, j \leq l, \quad 1 \leq \mu, \nu \leq k.
\end{equation*}

\noindent When describing an object with unknown grading indexed by $a$ the
weight will be described generically as $\varepsilon_a$.  This should
not be assumed to be an even weight.

\

The even positive roots of $osp(m|n)$ are composed entirely of the
usual positive roots of $o(m)$ together with those of $sp(n)$, namely:

\begin{alignat*}{3}
&\varepsilon_i \pm \varepsilon_j, & \qquad & 1 \leq i < j \leq l, \\
&\varepsilon_i, && 1 \leq i \leq l &&\quad \text{when } m=2l+1, \\
&\delta_\mu + \delta_\nu, && 1 \leq \mu,\nu \leq k, \\ &\delta_\mu -
\delta_\nu, && 1 \leq \mu < \nu \leq k.
\end{alignat*}

\noindent The root system also contains a set of odd positive roots,
which are:

\begin{equation*}
\delta_\mu + \varepsilon_i, \qquad 1 \leq \mu \leq k,\;1 \leq i \leq
m.
\hspace{2cm}
\end{equation*}

\noindent Throughout we choose to use the following set of simple
roots:

\begin{align*}
&\alpha_i = \varepsilon_i - \varepsilon_{i+1}, \hspace{11mm} 1 \leq i
  < l, \notag \\ &\alpha_l =
\begin{cases} 
\varepsilon_{l} + \varepsilon_{l-1},\quad & m=2l, \\ \varepsilon_l,
&m=2l+1,
\end{cases} \notag \\
&\alpha_\mu = \delta_\mu - \delta_{\mu+1}, \hspace{9mm} 1 \leq \mu <
k,\notag\\ &\alpha_s = \delta_k - \varepsilon_1.
\end{align*}

\noindent Note this choice is only valid for $m >2$.

Corresponding to these simple roots we have raising generators $e_a$,
lowering generators $f_a$ and Cartan elements $h_a$ given by:

\begin{alignat*}{4}
&e_i = \sigma^i_{i+1}, &\quad& f_i = \sigma^{i+1}_i, &\quad& h_i =
\sigma^i_i - \sigma^{i+1}_{i+1}, &\quad&1 \leq i < l, \notag \\ &e_l =
\sigma^{l-1}_{\overline{l}}, && f_l = \sigma^{\overline{l}}_{l-1}, &&
h_l = \sigma^{l-1}_{l-1} + \sigma^l_l, &&m=2l, \notag \\ &e_l =
\sigma^l_{l+1}, && f_l = \sigma^{l+1}_l, && h_l = \sigma^l_l, &&
m=2l+1, \notag \\ &e_\mu = \sigma^\mu_{\mu+1}, && f_\mu =
-\sigma^{\mu+1}_\mu, && h_\mu = \sigma^{\mu+1}_{\mu+1} -
\sigma^\mu_{\mu}, && 1 \leq \mu < k, \notag \\ &e_s =
\sigma^{\mu=k}_{i=1}, &&f_s = -\sigma^{i=1}_{\mu=k}, &&h_s =
-\sigma^{\mu=k}_{\mu=k} - \sigma^{i=1}_{i=1}.&&
\end{alignat*}

\noindent These automatically satisfy the defining relations of a Lie
superalgebra, which for the case at hand are:

\begin{alignat}{2} 
&[h_a, e_b] = (\alpha_a, \alpha_b) e_b, && \notag \\ &[h_a, f_b] = -
(\alpha_a, \alpha_b) f_b, && \notag \\ &[h_a, h_b] = 0,&& \notag \\
&[e_a, f_b] = \delta^a_b h_a, \notag \\ &[e_a, e_a] = [f_a,f_a]=0 &
\quad &\text{for } (\alpha_a, \alpha_a)=0,\notag\\ &(ad\,e_b\,
\circ)^{1-a_{bc}} e_c = 0 && \text{for } b \neq c,\; (\alpha_b,\alpha_b) 
  \neq 0, \label{Ser1} \\
&(ad\,f_b\, \circ)^{1-a_{bc}} f_c = 0 && \text{for } b \neq c, \;
  (\alpha_b,\alpha_b) \neq 0, \label{Ser2}
\end{alignat}
where the $a_{bc}$ are the entries of the corresponding Cartan matrix,

\begin{equation*}
a_{bc} =
\begin{cases}
   \frac{2(\alpha_b, \alpha_c)}{(\alpha_b, \alpha_b)}, \quad &
(\alpha_b, \alpha_b) \neq 0, \\ (\alpha_b, \alpha_c), &(\alpha_b,
\alpha_b) = 0,
\end{cases}
\end{equation*}

\noindent and $ad$ represents the \textit{adjoint action}

\begin{equation} \label{ad}
ad\, x \circ y = [x,y].
\end{equation}

The relations \eqref{Ser1} and \eqref{Ser2} are known as the
\textit{Serre relations} \cite{Serre}.  Superalgebras also have higher
order defining relations, not included here, which are known as the
\textit{extra Serre relations}.  They are dependent on the structure
of the chosen simple root system \cite{Yamane}.


\section{The $q$-Deformation: $U_q[osp(m|n)]$} \label{q-def}  

\noindent A quantum superalgebra is a generalised version of a
classical superalgebra involving a complex parameter $q$, which
reduces to the classical case as $q \rightarrow 1$.  In particular, we
construct $U_q[osp(m|n)]$ by \textit{$q$-deforming} the original
enveloping algebra of $osp(m|n)$ so that the generators remain
unchanged, but are now related by a deformation of the defining
relations.  Throughout this Communication $q$ is assumed not to be a
root of unity.

First note that in the enveloping algebra of $osp(m|n)$ the graded commutator
is realised by

\begin{equation*}
[A,B] = AB - (-1)^{[A][B]} BA.
\end{equation*}

\noindent With this operation, we then have:

\begin{definition} \label{def} The defining relations for $U_q[osp(m|n)]$ are:

\begin{alignat}{2}
&[h_a, e_b] = (\alpha_a, \alpha_b) e_b, && \notag \\ &[h_a, f_b] = -
(\alpha_a, \alpha_b) f_b, && \notag \\ &[h_a, h_b] = 0,&& \notag \\
&[e_a, f_b] = \delta^a_b \frac{(q^{h_a} - q^{-h_a})}{(q - q^{-1})},&&
\notag \\ &[e_a, e_a] = [f_a,f_a]=0 & \quad &\text{for } (\alpha_a,
\alpha_a)=0,\notag\\ &(ad\,e_b\, \circ)^{1-a_{bc}} e_c = 0   &&
\text{for }b \neq c,\; (\alpha_b,\alpha_b) \neq 0, \label{qS1}\\ &(ad\,f_b\, \circ)^{1-a_{bc}} f_c =
0 && \text{for } b \neq c, \; (\alpha_b,\alpha_b) \neq 0 \label{qS2}.
\end{alignat}

\end{definition}

The relations \eqref{qS1} and \eqref{qS2} are called the
\textit{$q$-Serre relations}.  Again, there are also extra $q$-Serre
relations which are not included here.  A complete list of them,
including those for affine superalgebras, can be found in
\cite{Yamane}.  Both the standard and extra $q$-Serre relations depend
on the adjoint action, which is no longer simply the graded commutator.  To
define the adjoint action for a quantum superalgebra, we first need
some new operations.

The \textit{coproduct}, $\Delta: U_q[osp(m|n)]^{\otimes 2} \rightarrow
U_q[osp(m|n)]^{\otimes 2}$, is the superalgebra homomorphism given by:

\begin{align}
&\Delta (e_a) = q^{\frac{1}{2}h_a} \otimes e_a + e_a \otimes
  q^{-\frac{1}{2} h_a}, \notag \\ &\Delta (f_a) = q^{\frac{1}{2}h_a}
  \otimes f_a + f_a \otimes q^{-\frac{1}{2} h_a},\notag \\ &\Delta
  (q^{\pm \frac{1}{2}h_a}) = q^{\pm \frac{1}{2}h_a} \otimes q^{\pm
  \frac{1}{2}h_a}, \notag \\ &\Delta (ab) = \Delta(a)
  \Delta(b). \label{coprod}
\end{align}

\noindent Note that in a $\mathbb{Z}_2$-graded algebra, multiplying tensor
products induces a grading term, according to

\begin{equation*}
(a \otimes b) (c \otimes d) = (-1)^{[b][c]} (ac \otimes bd).
\end{equation*}

\noindent We also require the \textit{antipode}, $S: U_q[osp(m|n)] \rightarrow
U_q[osp(m|n)]$, a superalgebra anti-homomorphism defined by:

\begin{align*}
&S (e_a) = - q^{-\frac{1}{2}(\alpha_a, \alpha_a)} e_a, \notag \\ &S
(f_a) = - q^{\frac{1}{2}(\alpha_a, \alpha_a)} f_a, \notag \\ &S
(q^{\pm h_a}) = q^{\mp h_a}, \notag \\ &S (ab) = (-1)^{[a][b]} S(b)
S(a).
\end{align*}

It can be shown that both the coproduct and antipode are consistent
with the defining relations of the superalgebra.  These mappings are
necessary to define the adjoint action for a quantum superalgebra, as
it can no longer be written simply in terms of the graded commutator.  If we
adopt Sweedler's notation for the coproduct \cite{Sweedler},

\begin{equation*}
\Delta(a) = \sum_{(a)} a^{(1)} \otimes a^{(2)},
\end{equation*}

\noindent the \textit{adjoint action} of $a$ on $b$ is defined to be

\begin{equation} \label{adj}
ad\; a \circ b = \sum_{(a)} (-1)^{[b][a^{(2)}]} a^{(1)} b S(a^{(2)}).
\end{equation}

\noindent in the case where $a$ is a simple generator.  Note that as $q \rightarrow 1$, equation \eqref{adj} reduces
to \eqref{ad}.

\

One quantity that repeatedly arises in calculations for both classical
and quantum Lie superalgebras is $\rho$, the \textit{graded half-sum
of positive roots}.  In the case of $U_q[osp(m|n)]$ it is given by:

\begin{equation*}
\rho = \frac{1}{2} \sum_{i=1}^l (m-2i) \varepsilon_i + \frac{1}{2}
\sum_{\mu=1}^k (n-m+2-2\mu) \delta_\mu.
\end{equation*}
This satisfies the property $(\rho, \alpha) = \frac{1}{2} (\alpha,
\alpha)$ for all simple roots $\alpha$.

\

As mentioned earlier, this root system and set of generators is only
valid for $m > 2, n \geq 0$.  When $m=0$, $U_q[osp(m|n)]$ is isomorphic to
$U_q[sp(n)]$.  Similarly, in \cite{Zhang} it was shown that every
finite-dimensional representation of $U_q[osp(1|n)]$ is isomorphic to
a finite-dimensional representation of $U_{-q}[o(n+1)]$.  As we are
only interested in finite-dimensional representations, and the
representation theory of these non-super quantum groups is
well-understood, we need not consider the cases with $m<2$.  Thus
although our root system is only valid for $m >2$, finding the Lax
operator for this root system will actually complete the work for all
$B$ and $D$ type quantum superalgebras.  This has, of course, already
been done for the more straightforward $A$ type quantum supergroups,
$U_q[gl(m|n)]$ \cite{Zhang2}.  The Lax operator for the $C$ type
quantum superalgebras ($U_q[osp(2|n)]$ ) is given in \cite{CLax}.

\subsection{$U_q[osp(m|n)]$ as a Quasi-Triangular Hopf Superalgebra}

\noindent A quantum superalgebra is actually a specific type of
quasi-triangular Hopf superalgebra.  This guarantees the existence of
a universal $R$-matrix, which provides a solution to the quantum
Yang--Baxter equation.  Before elaborating, we need to introduce the
graded twist map.

 The \textit{graded twist map} $T:U_q[osp(m|n)]^{\otimes 2}
 \rightarrow U_q[osp(m|n)]^{\otimes 2}$ is given by

\begin{equation*}
T(a \otimes b) = (-1)^{[a][b]} (b \otimes a).
\end{equation*}

\noindent For convenience $T \circ \Delta$, the twist map composed with the
coproduct, is denoted $\Delta^T$.  Then a \textit{universal
$R$-matrix}, $\mathcal{R}$, is an even, non-singular element of
$U_q[osp(m|n)]^{\otimes 2}$ satisfying the following properties:

\begin{align}
&\mathcal{R} \Delta (a) = \Delta^T (a)\mathcal{R}, \quad \forall a \in
  U_q[osp(m|n)], \notag \\ &(\text{id} \otimes \Delta) \mathcal{R} =
  \mathcal{R}_{13} \mathcal{R}_{12}, \notag \\ &(\Delta \otimes
  \text{id}) \mathcal{R} = \mathcal{R}_{13} \mathcal{R}_{23}.
  \label{Requations}
\end{align}

\noindent Here $\mathcal{R}_{ab}$ represents a copy of $\mathcal{R}$
acting on the $a$ and $b$ components respectively of $U_1 \otimes U_2
\otimes U_3$, where each $U$ is a copy of the quantum superalgebra
$U_q[osp(m|n)]$.  When $a>b$ the usual grading term from the twist map
is included, so for example $\mathcal{R}_{21} = [\mathcal{R}^T]_{12}
$, where $\mathcal{R}^T = T (\mathcal{R})$ is the \textit{opposite
universal $R$-matrix}.

One of the reasons $R$-matrices are significant is that as a
consequence of \eqref{Requations} they satisfy the YBE, which is prominent 
in the study of integrable systems \cite{Baxter}:

\begin{equation*}
\mathcal{R}_{12} \mathcal{R}_{13} \mathcal{R}_{23} = \mathcal{R}_{23}
  \mathcal{R}_{13} \mathcal{R}_{12}
\end{equation*}

A superalgebra may contain many different universal $R$-matrices, but
there is always a unique one belonging to $U_q[osp(m|n)]^- \otimes
U_q[osp(m|n)]^+$, and its opposite $R$-matrix in $U_q[osp(m|n)]^+
\otimes U_q[osp(m|n)]^-$.  Here $U_q[osp(m|n)]^-$ is the Hopf
subsuperalgebra generated by the lowering generators and Cartan
elements, while $U_q[osp(m|n)]^+$ is generated by the raising
generators and Cartan elements.  These particular $R$-matrices arise
out of Drinfeld's double construction \cite{Drinfeld}.  In this
Communication we consider the universal $R$-matrix belonging to
$U_q[osp(m|n)]^- \otimes U_q[osp(m|n)]^+$.

\section{An Ansatz for the Lax Operator}

\noindent Now we have the necessary information for the construction of a
\textit{Lax operator} for $U_q[osp(m|n)]$.  Previously this had only
been done in the superalgebra case for $U_q[gl(m|n)]$
\cite{Zhang}. Before defining a Lax operator, however, we need to
introduce the vector representation.

Let $\text{End} \; V$ be the set of endomorphisms of $V$, an
$(m+n)$-dimensional vector space.  Then the irreducible \textit{vector
representation} $\pi: U_q[osp(m|n)] \rightarrow \text{End} \; V$ is
left undeformed from the classical vector representation of
$osp(m|n)$, which acts on the Cartan-Weyl generators given in equation
\eqref{CW} according to:

\begin{equation*}
\pi (\sigma^a_b) = E^a_b - (-1)^{[a]([a]+[b])} \xi_a \xi_b
E^{\overline{b}}_{\overline{a}}
\end{equation*}

\noindent where $E^a_b$ is the $(m+n) \times (m+n)$-dimensional
elementary matrix with $(a,b)$ entry $1$ and zeroes elsewhere.

Now let $\mathcal{R}$ be a universal $R$-matrix of $U_{q}[osp(m|n)]$
and $\pi$ the vector representation.  The Lax operator associated with
$\mathcal{R}$ is given by

\begin{equation*}
R = (\pi \otimes \text{id}) \mathcal{R} \in (\text{End} \; V) \otimes
U_{q}[osp(m|n)]
\end{equation*}

\noindent and the $R$-matrix in the vector representation $\mathfrak{R}$
is given by:

\begin{equation*}
\mathfrak{R} = (\pi \otimes \pi) \mathcal{R} = (\text{id} \otimes \pi) R
\in (\text{End} \; V) \otimes (\text{End} \; V).
\end{equation*}

\noindent Then the Yang-Baxter equation reduces to:

\begin{equation*}
\mathfrak{R}_{12} R_{13} R_{23} = R_{23} R_{13} \mathfrak{R}_{12}
\end{equation*}

\noindent acting on the space $V \otimes V \otimes U_q[osp(m|n)]$.

Previously an $R$-matrix in the vector representation, $\mathfrak{R}$, has
been calculated for both $U_q[osp(m|n)]$ and its affine extension
\cite{Scheun,Martins,Mehta}, however in general the Lax operator is still
unknown.  The Lax operator is significant because we can use it to
calculate solutions to the quantum Yang--Baxter equation for an
arbitrary finite-dimensional representation.

In the following sections we also sometimes make use of the bra and
ket notation.  The set $\{| a \rangle, a=1,...,m+n \}$ is a basis for
$V$ satisfying the property

\begin{equation*}
E^a_b |c \rangle = \delta^c_b |a \rangle.
\end{equation*}

\noindent The set $\{ \langle a|,a=1,...,m+n \}$ is the dual basis
such that

\begin{equation*}
\langle c| E^a_b = \delta^a_c \langle b| \quad \text{and} \quad
\langle a|b \rangle = \delta^a_b.
\end{equation*}   

\

\noindent As we wish to find the Lax operator belonging to $ \pi \bigl(
U_q[osp(m|n)]^-\bigr) \otimes U_q[osp(m|n)]^+$, we adopt the following
ansatz for $R$:

\begin{equation*} \label{RmatrixAnsatz}
R \equiv q ^ {\underset{a}{\sum} \pi (h_{a}) \otimes h^{a}} \Bigl[ I \otimes
I + (q - q^{-1}) \sum_{\varepsilon_{a} < \varepsilon_{b}} (-1)^{[b]}
E^a_b \otimes \hat{\sigma}_{ba} \Bigr].
\end{equation*} 

\noindent Here $\{ h_a \}$ is a basis for the Cartan subalgebra such
that $h_a = h_{\varepsilon_a}$, and $\{ h^a \}$ the dual basis, so
$h^a = (-1)^{[a]} h_{\varepsilon_a}$. The $\hat{\sigma}_{ba}$ are
unknown operators for which we are trying to solve.  Throughout the remainder 
of this Communication, when working in the vector representation, we simply
use $h_a$ rather than $\pi(h_a)$, and $e_a$ rather than $\pi(e_a)$.
\subsection{Constraints Arising from the Defining Relations} \label{Developing Relations}

\noindent The Lax operator $R$ must be consistent with the defining
relations for the $R$-matrix, which were given as equation
\eqref{Requations}.  In particular, it must satisfy the intertwining
property for the raising generators,

\begin{equation*}
R \Delta (e_c) = \Delta^T (e_c) R.
\end{equation*}

\noindent To apply this, recall that

\begin{equation*}
\Delta(e_c) = q^{\frac{1}{2} h_c} \otimes e_c + e_c \otimes
  q^{-\frac{1}{2} h_c}.
\end{equation*}

\noindent Then, from the defining relations, we have

\begin{align*}
\Delta^T(e_c) q^{\underset{a}{\sum} h_a \otimes h^a} &=
  (q^{\frac{1}{2} h_c} \otimes e_c + e_c \otimes q^{-\frac{1}{2} h_c})
  q^{\underset{a}{\sum} h_a \otimes h^a} \\ &= q^{\underset{a}{\sum}
  h_a \otimes h^a} (e_c \otimes q^{-\frac{1}{2} h_c} + q^{-\frac{3}{2}
  h_c} \otimes e_c).
\end{align*}

\noindent Using this, we see

\begin{align}
\Delta^T(e_c)R &= q^{\underset{a}{\sum} h_a \otimes h^a} (e_c \otimes
  q^{-\frac {1}{2} h_c} + q^{-\frac{3}{2} h_c} \otimes e_c) \notag \\
  &\hspace{4cm} \times \Bigl[ I \otimes I + (q - q^{-1})
  \sum_{\varepsilon_{a}< \varepsilon_{b}} (-1)^{[b]} E^a_b \otimes
  \hat{\sigma}_{ba} \Bigr] \notag \\ &= q^{\underset{a}{\sum} h_a
  \otimes h^a} \biggl\{ e_c \otimes q^{-\frac{1}{2} h_c} +
  q^{-\frac{3}{2} h_c} \otimes e_c \notag \\ &\hspace{22mm} + (q -
  q^{-1}) \sum_{\varepsilon_{a} < \varepsilon_{b}} (-1)^{[b]} \Bigl[
  e_c E^a_b \otimes q^{-\frac{1}{2} h_c} \hat{\sigma}_{ba} \notag \\ &
  \hspace{45mm} + (-1)^{([a] + [b])[c]} q^{-\frac{3}{2} (\alpha_c,
  \varepsilon_a)} E^a_b \otimes e_c \hat{\sigma}_{ba} \Bigr] \biggr\}.
   \label{delR} 
\end{align}

\noindent Also,

\begin{align}
R \Delta (e_c) &= q^{\underset{a}{\sum}h_a \otimes h^a} \biggl\{
   q^{\frac{1}{2} h_c} \otimes e_c + e_c \otimes q^{-\frac{1}{2} h_c}
   \notag \\ & \hspace{21mm} + (q - q^{-1}) \sum_{\varepsilon_{a} <
   \varepsilon_{b}} (-1)^{[b]} \Bigl[ q^{\frac{1}{2} (\alpha_c,
   \varepsilon_b)} E^a_b \otimes \hat{\sigma}_{ba} e_c \notag \\ &
   \hspace{50mm} + (-1)^{([a]+[b])[c]} E^a_b e_c \otimes
   \hat{\sigma}_{ba} q^{-\frac{1}{2} h_c} \Bigr] \biggr\}
   \label{Rdel}.
\end{align}


Hence to apply the intertwining property we simply equate (\ref{delR})
and (\ref{Rdel}).  First note that $R$ is weightless, so
$\hat{\sigma}_{ba}$ has weight $\varepsilon_b - \varepsilon_{a}$, and
thus

\begin{equation*}
q^{-\frac{1}{2} h_c} \hat{\sigma}_{ba} = q^{-\frac{1}{2} (\alpha_c,
  \varepsilon_b - \varepsilon_a)} \hat{\sigma}_{ba} q^{-\frac{1}{2}
  h_c}.
\end{equation*}

\noindent Then, equating those terms with zero weight in the first
element of the tensor product, we obtain

\begin{multline}
(q^{\frac{1}{2} h_c} - q^{-\frac{3}{2} h_c}) \otimes e_c \\ =
  (q-q^{-1}) \sum_{\varepsilon_{b} - \varepsilon_{a} = \alpha_c}
  (-1)^{[b]} \bigl( q^{-\frac{1}{2} (\alpha_c, \alpha_c)} e_c E^a_b -
  (-1)^{[c]} E^a_b e_c \bigr) \otimes \hat{\sigma}_{ba}
  q^{-\frac{1}{2} h_c}.\label{eq1}
\end{multline}

\noindent Comparing the remaining terms, we also find

\begin{multline} \label{eq2}
\underset{\varepsilon_{b} - \varepsilon_{a} \neq \alpha_c}
  {\sum_{\varepsilon_ {a} < \varepsilon_{b}}} (-1)^{[b]} \bigl(
  q^{-\frac{1}{2}(\alpha_c, \varepsilon_{b} - \varepsilon_{a})} e_c
  E^a_b - (-1)^{([a] + [b])([c])} E^a_b e_c \bigr) \otimes
  \hat{\sigma}_{ba} q^{-\frac{1}{2} h_c} \\ = \sum_{\varepsilon_{a} <
  \varepsilon_{b}} (-1)^{[b]} E^a_b \otimes \bigl( q^{\frac{1}{2}
  (\alpha_c, \varepsilon_b)} \hat{\sigma}_{ba} e_c -
  (-1)^{([a]+[b])[c]} q^{-\frac{3}{2} (\alpha_c, \varepsilon_a)} e_c
  \hat{\sigma}_{ba} \bigr).
\end{multline}



\

 From the first of these equations we can deduce the simple values of
 $\hat{\sigma}_{ba}$, namely those
for which $\varepsilon_b - \varepsilon_a$ is a simple root; 
 from the second, relations involving all the
 $\hat{\sigma}_{ba}$.  Before doing so, however, it is convenient to
 define a new set, $\overline{\Phi}^+$.

\begin{definition}
The extended system of positive roots, $\overline{\Phi}^+$, is defined
by

\begin{equation*}
\overline{\Phi}^+ \equiv \{\varepsilon_b - \varepsilon_a|
\varepsilon_b > \varepsilon_a\} = \Phi^+ \cup \{2 \varepsilon_i | 1
\leq i \leq l\}
\end{equation*}

\noindent where $\Phi^+$ is the usual system of positive roots.
\end{definition}

\

Now consider equation (\ref{eq2}).  In the case when $\varepsilon_b -
\varepsilon_a + \alpha_c \notin \overline{\Phi}^+$, by collecting the
terms of weight $\varepsilon_b - \varepsilon_a + \alpha_c$ in the
second half of the tensor product we find:

\begin{equation} \label{*}
q^{\frac{1}{2}(\alpha_c, \varepsilon_{b})} \hat{\sigma}_{ba} e_c -
    (-1)^{([a]+[b])[c]} q^{-\frac{3}{2} (\alpha_c, \varepsilon_a)} e_c
    \hat{\sigma}_{ba} = 0.
\end{equation}

\noindent Similarly, when $\varepsilon_b > \varepsilon_a$ and
$\varepsilon_b - \varepsilon_a + \alpha_c = \varepsilon_{b'} -
\varepsilon_{a'} \in \overline{\Phi}^+$ we find:

\begin{multline*}
\underset{\varepsilon_{b} - \varepsilon_{a} + \alpha_c =
  \varepsilon_{b'} - \varepsilon_{a'}} {\sum_{\varepsilon_{a'} <
  \varepsilon_{b'}}} \hspace{-8mm} (-1)^{[b']}
  \bigl(q^{-\frac{1}{2}(\alpha_c, \varepsilon_{b'}-\varepsilon_{a'})}
  e_c E^{a'}_{b'} - (-1)^{([a']+[b'])[c]} E^{a'}_{b'} e_c \bigr)
  \otimes \hat{\sigma}_{b'\!a'} q^{-\frac{1}{2} h_c} \\ = (-1)^{[b]}
  E^a_b \otimes \bigl( q^{\frac{1}{2} (\alpha_c, \varepsilon_b)}
  \hat{\sigma}_{ba} e_c - (-1)^{([a]+[b])[c])} q^{-\frac{3}{2}
  (\alpha_c, \varepsilon_a)} e_c \hat{\sigma}_{ba} \bigr).
\end{multline*}

\noindent However $e_c E^{a'}_{b'}$ and $E^a_{b}$ are linearly
independent unless $b=b'$, as are $E^{a'}_{b'} e_c$ and $E^{a}_{b}$
for $a \neq a'$, and thus this equation reduces to

\begin{multline*}
q^{-\frac{1}{2} (\alpha_c, \varepsilon_b - \varepsilon_{a} +
   \alpha_c)} e_c E^{a'}_{b} \otimes \hat{\sigma}_{ba'}
   q^{-\frac{1}{2} h_c} \Big\vert_{\varepsilon_{a'} = \varepsilon_a -
   \alpha_c} \\
 \hspace{-12mm}- (-1)^{([a]+[b])[c]} E^{a}_{b'}e_c \otimes
  \hat{\sigma}_{b'a} q^{-\frac{1}{2}h_c} \Big\vert_{\varepsilon_{b'} =
  \varepsilon_b +\alpha_c}\\ = \; E^{a}_{b} \otimes \bigl(
  q^{\frac{1}{2} (\alpha_c,\varepsilon_b)} \hat{\sigma}_{ba} e_c -
  (-1)^{([a]+[b])[c]} q^{-\frac{3}{2} (\alpha_c,\varepsilon_a)} e_c
  \hat{\sigma}_{ba} \bigr)
\end{multline*}

\noindent for $\varepsilon_b > \varepsilon_a$.  This further
simplifies to

\begin{multline}
q^{-\frac{1}{2} (\alpha_c, \alpha_c - \varepsilon_{a})} \langle
  a|e_c|a' \rangle \hat{\sigma}_{ba'}-(-1)^{([a]+[b])[c]}
  q^{\frac{1}{2} (\alpha_c, \varepsilon_b)} \langle b'|e_c|b \rangle
  \hat{\sigma}_{b'a} \\ = q^{(\alpha_c,\varepsilon_b)}
  \hat{\sigma}_{ba}e_c q^{\frac{1}{2} h_c} - (-1)^ {([a]+[b])[c]}
  q^{-(\alpha_c,\varepsilon_a)} e_c q^{\frac{1}{2} h_c}
  \hat{\sigma}_{ba} \label{**}
\end{multline}

\noindent for $\varepsilon_b > \varepsilon_a$.  All the necessary
information is contained within equations \eqref{eq1} and \eqref{**}.
To construct the Lax operator $R = (\pi \otimes 1) \mathcal{R}$ first
we use equation \eqref{eq1} to find the solutions for
$\hat{\sigma}_{ba}$ associated with the simple roots $\alpha_c$.  Then
we apply the recursion relations arising from \eqref{**} to find the
remaining values of $\hat{\sigma}_{ba}$.


\subsection{The Simple Operators} \label{iv}

\noindent In this section we solve equation \eqref{eq1}, rewritten
below, to find the simple values of $\hat{\sigma}_{ba}$.

\begin{multline*} 
(q^{\frac{1}{2} h_c} - q^{-\frac{3}{2} h_c}) \otimes e_c \\ =
(q-q^{-1}) \sum_{\varepsilon_{b} - \varepsilon_{a} = \alpha_c}
(-1)^{[b]} \bigl( q^{-\frac{1}{2} (\alpha_c, \alpha_c)} e_c E^a_b -
(-1)^{[c]} E^a_b e_c \bigr) \otimes \hat{\sigma}_{ba} q^{-\frac{1}{2}
h_c}. \tag{\ref{eq1}}
\end{multline*}

\

To solve this we must consider the various simple roots individually.
Consider the simple roots $\alpha_i = \varepsilon_i -
\varepsilon_{i+1},\; 1 \leq i < l$.  In the vector representation we
know
$$e_i = E^i_{i+1} - E^{\overline{i+1}}_{\overline{i}}, \quad h_i =
E^i_i - E^{i+1}_{i+1} + E^{\overline{i+1}}_{\overline{i+1}} -
E^{\overline{i}}_{\overline{i}}.$$

\noindent Hence the left-hand side of (\ref{eq1}) becomes:

\begin{align*}
LHS &= (q^{\frac{1}{2} h_i} - q^{-\frac{3}{2} h_i}) \otimes e_i \notag
    \\ &= \bigl\{ (q^{\frac{1}{2}} - q^{-\frac{3}{2}})(E^i_i +
    E^{\overline{i+1}} _{\overline{i+1}}) + (q^{-\frac{1}{2}} -
    q^{\frac{3}{2}})(E^{i+1}_{i+1} + E^{\overline{i}}_{\overline{i}})
    \bigr\} \otimes e_i \notag \\ &= (q-q^{-1}) \bigl\{
    q^{-\frac{1}{2}} (E^i_i + E^{\overline{i+1}}_ {\overline{i+1}}) -
    q^{\frac{1}{2}} (E^{i+1}_{i+1} + E^{\overline{i}} _{\overline{i}})
    \bigr\} \otimes e_i,
\end{align*}

\noindent whereas the right-hand side is:

\begin{align*}
RHS &= (q-q^{-1}) \sum_{\varepsilon_{b} - \varepsilon_{a} = \alpha_i}
   \bigl( q^{-1} e_i E^a_b - E^a_b e_i \bigr) \otimes
   \hat{\sigma}_{ba} q^{-\frac{1}{2} h_i} \notag \\ &= (q-q^{-1})
   \Bigl\{ (q^{-1} E^i_i - E^{i+1}_{i+1}) \otimes \hat{\sigma}_{i\,
   i+1} q^{-\frac{1}{2} h_i}\\
& \hspace{55mm} - (q^{-1}  E^{\overline{i+1}}_{\overline{i+1}} -
   E^{\overline{i}}_{\overline{i}}) \otimes
   \hat{\sigma}_{\overline{i+1}\, \overline{i}} q^{-\frac{1}{2} h_i}
   \Bigr\}.
\end{align*}

\noindent Equating these gives

\begin{equation*}
\hat{\sigma}_{i\, i+1} = - \hat{\sigma}_{\overline{i+1}\,
 \overline{i}} = q^{\frac{1}{2}} e_i q^{\frac{1}{2} h_i}, \quad 1 \leq
 i < l.
\end{equation*}

By performing similar calculations for the other simple roots we
obtain the simple operators given below in Table \ref{fundval}.  These
values for $\hat{\sigma}_{ba}$ form the basis for finding $R$, as from
these all the others can be explicitly determined in any given
representation.

\begin{table}[ht]
\caption{The simple values for $\hat{\sigma}_{ba}$.} 
\label{fundval}
\centering
\begin{tabular}{|l|lll|}
\hline\noalign{\smallskip}
\multicolumn{1}{|c|}{Simple Root}& \multicolumn{3}{c|}{Corresponding
  $\hat{\sigma}_{ba}$} \\ 
\noalign{\smallskip}\hline\noalign{\smallskip}
$\alpha_i = \varepsilon_i -
  \varepsilon_{i+1},\,i < l$ & $\hat{\sigma}_{i\, i+1}$&$ =
  -\hat{\sigma}_{\overline{i+1}\,\overline{i}}$&$ = q^{\frac{1}{2}}
  e_i q^{\frac{1}{2} h_i}$ \\ $\alpha_l = \varepsilon_{l-1} +
  \varepsilon_l,\,m = 2l$ & $\hat{\sigma}_ {l-1\,\overline{l}}$&$ =
  -\hat{\sigma}_{l\, \overline{l-1}}$&$ = q^{\frac{1}{2}} e_l
  q^{\frac{1}{2} h_l}$ \\ $\alpha_l = \varepsilon_l,\,m=2l+1$ &
  $\hat{\sigma}_{l\, l+1}$&$ = -q^{-\frac {1}{2}} \hat{\sigma}_{l+1\,
  \overline{l}}$&$ =e_l q^{\frac{1}{2} h_l} $ \\ $\alpha_\mu =
  \delta_\mu - \delta_{\mu+1},\,\mu < k$ & $\hat{\sigma}_ {\mu\,
  \mu+1} $&$ = \hat{\sigma}_{\overline{\mu+1}\, \overline{\mu}} $&$=
  q^{-\frac{1}{2}} e_\mu q^{\frac{1}{2} h_\mu}$ \\ $\alpha_s =
  \delta_k - \varepsilon_1,$ & $\hat{\sigma}_{\mu=k\, i=1} $&$ =
  (-1)^k q \hat{\sigma}_{\overline{i}=\overline{1}\, \overline{\mu} =
  \overline{k}} $&$ = q^{\frac{1}{2}} e_s q^{\frac{1}{2} h_s}$\\
\noalign{\smallskip}\hline
\end{tabular}
\end{table}

\subsection{Constructing the Non-Simple Operators} \label{noniv}

\noindent Now we develop the recurrence relations required to
calculate the remaining values of $\hat{\sigma}_{ba}$.  Recall that
for \mbox{$\varepsilon_b > \varepsilon_a$},

\begin{multline} \label{**2}
q^{-\frac{1}{2} (\alpha_c, \alpha_c - \varepsilon_{a})} \langle
  a|e_c|a' \rangle \hat{\sigma}_{ba'}-(-1)^{([a]+[b])[c]}
  q^{\frac{1}{2} (\alpha_c, \varepsilon_b)} \langle b'|e_c|b \rangle
  \hat{\sigma}_{b'a} \\ = q^{(\alpha_c,\varepsilon_b)}
  \hat{\sigma}_{ba} e_c q^{\frac{1}{2} h_c} - (-1)^ {([a]+[b])[c]}
  q^{-(\alpha_c,\varepsilon_a)} e_c q^{\frac{1}{2} h_c}
  \hat{\sigma}_{ba}.
\end{multline}

To extract the recurrence relations to be applied to the simple values
of $\hat{\sigma}_{ba}$, we must again consider the simple roots
individually.  We begin with the case $\alpha_i = \varepsilon_i -
\varepsilon_{i+1}$, so \hbox{$e_i = \sigma^i_{i+1} \equiv E^i_{i+1} -
E^{\overline{i+1}}_{\overline{i}}$}.  Now

\begin{equation*}
\langle a|e_i = \delta_{ai} \langle i+1| - \delta_{a\, \overline{i+1}}
  \langle \overline{i}|,\qquad e_i |b \rangle = \delta_{b\, i+1} |i
  \rangle - \delta_{b \overline{i}} |\overline{i+1} \rangle.
\end{equation*}

\noindent We then apply this to equation \eqref{**2} to obtain:

\begin{align*}
q^{-\frac{1}{2}(\alpha_i, \alpha_i)} \bigl\{ \delta_{a i}
   q^{\frac{1}{2}(\alpha_i, \varepsilon_i)} \hat{\sigma}_{b\, i+1}-
   \delta_{a\, \overline{i+1}} &q^{-\frac{1}{2}(\alpha_i,
   \varepsilon_{i+1})} \hat{\sigma}_ {b\, \overline{i}} \bigr\} \notag
   \\ - \bigl\{ \delta_{b\, i+1} q^{\frac{1}{2} (\alpha_i,
   \varepsilon_{i+1})}&\hat{\sigma}_{ia} - \delta_{b\, \overline{i}}
   q^{-\frac{1}{2} (\alpha_i, \varepsilon_i)}
   \hat{\sigma}_{\overline{i+1}\,a} \bigr\} \notag \\ =
   &q^{(\alpha_i,\varepsilon_b)} \hat{\sigma}_{ba}e_i q^{\frac{1}{2}
   h_i} - q^{-(\alpha_i,\varepsilon_a) } e_i q^{\frac{1}{2} h_i}
   \hat{\sigma}_{ba}
\end{align*}

\noindent for $\varepsilon_b >\varepsilon_a.$ Recalling that
$\hat{\sigma}_{i\, i+1} = - \hat{\sigma}_{\overline{i+1}\,
\overline{i}} = q^{\frac{1}{2}} e_i q^{\frac{1}{2} h_i},\; 1 \leq i <
l$, the above simplifies to

\begin{align*}
\delta_{ai} \hat{\sigma}_{b\, i+1}\! -\delta_{a\, \overline{i+1}}
  \hat{\sigma}_ {b\overline{i}} -\delta_{b\,i+1} \hat{\sigma}_{ia}\!
  +\delta_ {b\overline{i}} \hat{\sigma}_{\overline{i+1}\,a} \!
  &=q^{(\alpha_i,\varepsilon_b)} \hat{\sigma}_{ba} \hat{\sigma}_{i\,
  i+1}\! - q^{-(\alpha_i,\varepsilon_a)} \hat{\sigma}_{i\, i+1}
  \hat{\sigma}_{ba} \\ &=q^{-(\alpha_i,\varepsilon_a)}
  \hat{\sigma}_{\overline{i+1}\, \overline{i}} \hat{\sigma}_{ba}\! -
  q^{(\alpha_i,\varepsilon_b)} \hat{\sigma}_{ba}
  \hat{\sigma}_{\overline{i+1}\, \overline{i}}.
\end{align*}

\noindent From this we can deduce the following relations for $1 \leq
i <l$:

\begin{alignat}{2}
&\hat{\sigma}_{b\, i+1} = \hat{\sigma}_{b\,i} \hat{\sigma}_{i\, i+1} -
   q^{-1} \hat{\sigma}_{i\, i+1} \hat{\sigma}_{b\, i}, & \qquad
   &\varepsilon_b > \varepsilon_i, \notag \\
   &\hat{\sigma}_{\overline{i+1}\, a} = \hat{\sigma}_{\overline{i+1}\,
   \overline{i}} \hat{\sigma}_{\overline{i}\, a} - q^{-1}
   \hat{\sigma}_ {\overline{i}\, a} \hat{\sigma}_{\overline{i+1}\,
   \overline{i}}, & & \varepsilon_a < -\varepsilon_i, \notag \\
   &\hat{\sigma}_{b\, \overline{i}} = q^{(\alpha_i, \varepsilon_b)}
   \hat{\sigma}_ {b\, \overline{i+1}} \hat{\sigma}_{\overline{i+1}\,
   \overline{i}} - q^{-1} \hat{\sigma}_{\overline{i+1}\,\overline{i}}
   \hat{\sigma}_{b\,\overline{i+1}} , & & \varepsilon_b > -
   \varepsilon_{i+1},\, b \neq i+1, \notag \\ &\hat{\sigma}_{i\,a} =
   q^{-(\alpha_i, \varepsilon_a)} \hat{\sigma}_{i\, i+1}
   \hat{\sigma}_{i+1\,a} - q^{-1} \hat{\sigma}_{i+1\,a}
   \hat{\sigma}_{i\,i+1}, & & \varepsilon_a < \varepsilon_{i+1}, \; a
   \neq \overline{i+1},\notag \\ &\hat{\sigma}_{i\, \overline{i+1}} +
   \hat{\sigma}_{i+1\,\overline{i}} = q^{-1} \bigl[
   \hat{\sigma}_{i\,i+1}, \hat{\sigma}_{i+1\,\overline{i+1}} \bigr],&&
   \notag \\ &q^{(\alpha_i, \varepsilon_b)} \hat{\sigma}_{ba}
   \hat{\sigma}_{i\,i+1} - q^{-(\alpha_i, \varepsilon_a)}
   \hat{\sigma}_{i\, i+1} \hat{\sigma}_{ba} = 0, && \varepsilon_b >
   \varepsilon_a,\; a \neq i, \overline{i+1} \notag \\ &&&
   \hspace{13mm} \text{ and }b \neq i+1, \overline{i}. \notag
\end{alignat}

We then follow the same procedure to find the relations associated
with the other simple roots.  A complete list of the relations derived
in this manner are given in Tables \ref{list}, \ref{even} and
\ref{odd} in Appendix A.

Although the list of relations is long, they can be summarised in
a compact form.  There are two different types of relations; recursive
and $q$-commutative.  The latter can be condensed into:

\vspace{-2mm}
\begin{equation} \label{commutationrelations}
\boxed{q^{(\alpha_c, \varepsilon_b)} \hat{\sigma}_{ba} e_c
  q^{\frac{1}{2} h_c} - (-1)^{([a]+[b])[c]} q^{-(\alpha_c,
    \varepsilon_a)} e_c q^{\frac{1}{2}h_c} \hat{\sigma}_{ba} = 0,
  \quad \varepsilon_b > \varepsilon_a}
\end{equation}
      
\noindent where neither $\varepsilon_a - \alpha_c$ nor $\varepsilon_b
+ \alpha_c$ equals any $\varepsilon_x$.  Note this is almost the same
as equation \eqref{*}, with slightly softer restrictions on $a$ and
$b$ the only difference.

The recursion relations can, in the case $m=2l+1$, be summarised as:

\vspace{-2mm}
\begin{equation} \label{nice}
\boxed{\hat{\sigma}_{ba} = q^{-(\varepsilon_b,\varepsilon_a)}
  \hat{\sigma}_{bc} \hat{\sigma}_{ca} - q^{-(\varepsilon_c,
  \varepsilon_c)} (-1)^{([b]+[c]) ([a]+[c])} \hat{\sigma}_{ca}
  \hat{\sigma}_{bc}, \quad \varepsilon_b > \varepsilon_c >
  \varepsilon_a}
\end{equation}

\noindent where $ c \neq \overline{b} \text{ or } \overline{a}$.

\noindent In the case $m=2l$, one extra relation is required, namely:

\begin{equation} \label{2l}
\hat{\sigma}_{l-1\, \overline{l}} + \hat{\sigma}_{l\,
  \overline{l-1}} = q^{-1} [ \hat{\sigma}_{l-1\, l}, \hat{\sigma}_{l
  \overline{l}} ],\quad m=2l.
\end{equation}

\

Noting from Table \ref{fundval} that $\hat{\sigma}_{l-1\,
\overline{l}} + \hat{\sigma}_{l\, \overline{l-1}} = 0$ and from the
Appendix that $\hat{\sigma}_{l\, \overline{l}}$ commutes with all the
other simple generators, equation \eqref{2l} can be replaced with the
condition:

\begin{equation} \label{l-1lbar}
\boxed{\hat{\sigma}_{l \overline{l}} =0,\quad m=2l.}
\end{equation}

\

It is not difficult to verify that all these recursion relations can
be obtained from those derived from equation \eqref{**2}.  To show the
reverse is tedious, but straightforward \cite{me}.  This unified form
of the relations can be applied to the simple operators given in Table
\ref{fundval} to obtain all the remaining values of
$\hat{\sigma}_{ba}$ in a given representation.

Hence we have found the following result:

\begin{proposition}

\noindent There is a unique operator-valued matrix $R \in (\text{End
}V) \otimes U_q[osp(m|n)]^+$ of the form

\begin{equation*}
R = q^{\underset{a}{\sum} h_{a} \otimes h^{a}} \Bigl[I \otimes I + (q
- q^{-1}) \sum_{\varepsilon_{a} < \varepsilon_{b}} (-1)^{[b]} E^a_b
\otimes \hat{\sigma}_{ba} \Bigr],
\end{equation*} 

\noindent satisfying \mbox{$R \Delta (e_c) = \Delta^T (e_c) R$}. The
simple operators for that matrix are given by:

\begin{alignat}{2}
&\hat{\sigma}_{i\, i+1} = - \hat{\sigma}_{\overline{i+1}\,
  \overline{i}} = q^{\frac{1}{2}} e_i q^{\frac{1}{2} h_i},&&1 \leq i <
  l, \notag \\ &\hat{\sigma}_{l-1\, \overline{l}} = -
  \hat{\sigma}_{l\, \overline{l-1}} = q^{\frac{1}{2}} e_l
  q^{\frac{1}{2} h_l}, && m=2l, \notag \\ &\hat{\sigma}_{l\,
  \overline{l}} = 0, && m=2l, \notag \\ &\hat{\sigma}_{l\, l+1} = -
  q^{-\frac{1}{2}} \hat{\sigma}_{l+1\, \overline{l}} = e_l
  q^{\frac{1}{2} h_l}, && m=2l+1, \notag \\ &\hat{\sigma}_{\mu \, \mu
  +1} = \hat{\sigma}_{\overline{\mu+1}\,\overline{\mu}} =
  q^{-\frac{1}{2}} e_\mu q^{\frac{1}{2} h_\mu}, &\qquad&1 \leq \mu <k,
  \notag \\ &\hat{\sigma}_{\mu=k \, i=1} = (-1)^k q\, \hat{\sigma}_{i
  = \overline{1} \, \overline{\mu} = \overline{k}} = q^{\frac{1}{2}}
  e_s q^{\frac{1}{2} h_s};
\label{ivalues}
\end{alignat}

\noindent and the remaining values can be calculated using

\

\noindent (i) the $q$-commutation relations

\begin{equation} \label{qcom}
q^{(\alpha_c, \varepsilon_b)} \hat{\sigma}_{ba} e_c q^{\frac{1}{2}
  h_c} - (-1)^{([a]+[b])[c]} q^{-(\alpha_c, \varepsilon_a)} e_c
  q^{\frac{1}{2}h_c} \hat{\sigma}_{ba} = 0, \quad \varepsilon_b >
  \varepsilon_a,
\end{equation}
      
\noindent where neither $\varepsilon_a - \alpha_c$ nor $\varepsilon_b
+ \alpha_c$ equals any $\varepsilon_x$; and

\

\noindent (ii) the induction relations

\begin{equation} \label{indrel}
\hat{\sigma}_{ba} = q^{-(\varepsilon_b,\varepsilon_a)}
  \hat{\sigma}_{bc} \hat{\sigma}_{ca} - q^{-(\varepsilon_c,
  \varepsilon_c)} (-1)^{([b]+[c]) ([a]+[c])} \hat{\sigma}_{ca}
  \hat{\sigma}_{bc}, \quad \varepsilon_b > \varepsilon_c >
  \varepsilon_a,
\end{equation}

\noindent where $ c \neq \overline{b} \text{ or } \overline{a}$.
\end{proposition}

\

We have found a set of simple operators and relations which uniquely
define the unknowns $\hat{\sigma}_{ba}$.  Theoretically the resultant
matrix $R$ must be a Lax operator, as we know there is one of the
given form.  We choose to confirm this, however, by verifying that $R$
satisfies the remaining $R$-matrix properties.  These are

\begin{equation} \label{1tensordel}
(\text{id} \otimes \Delta) R = R_{13} R_{12}
\end{equation}

\noindent and the intertwining property for the remaining generators,

\begin{equation*}
R \Delta(a) = \Delta^T(a) R, \qquad \forall a \in U_q[osp(m|n)].
\end{equation*}

\noindent We also calculate the opposite Lax operator $R^T$, and
briefly examine whether the defining relations for the
$\hat{\sigma}_{ba}$ incorporate the $q$-Serre relations for
$U_q[osp(m|n)]$.


\subsection{Calculating the Coproduct}
\noindent We begin by considering the first of these defining
properties, equation \eqref{1tensordel}.  In order to evaluate
$(\text{id} \otimes \Delta) R$, however, we need to know $\Delta
(\hat{\sigma}_{ba})$.  First note an alternative way of writing $R$ is

\begin{equation*}
R = \sum_a E^a_a \otimes q^{h_{\varepsilon_a}} + (q-q^{-1})
  \sum_{\varepsilon_a < \varepsilon_b} (-1)^{[b]} E^a_b \otimes
  q^{h_{\varepsilon_a}} \hat{\sigma}_{ba},
\end{equation*}

\noindent with the $\hat{\sigma}_{ba}$ as given before.  Using this
form, we find

\begin{align*}
R_{13} R_{12} &= \Bigl( \sum_a E^a_a \otimes I \otimes
    q^{h_{\varepsilon_a}} \; + \; (q-q^{-1}) \sum_{\varepsilon_{b} >
    \varepsilon_{a}} (-1)^{[b]} E^a_b \otimes I \otimes
    q^{h_{\varepsilon_a}} \hat{\sigma}_{ba} \Bigr) \notag \\ &
    \hspace{12mm} \Bigl( \sum_c E^c_c \otimes q^{h_{\varepsilon_c}}
    \otimes I \; + \; (q-q^{-1}) \sum_{\varepsilon_{d} >
    \varepsilon_{c}} (-1)^{[d]} E^c_d \otimes q^{h_{\varepsilon_c}}
    \hat{\sigma}_{dc} \otimes I \Bigr) \notag \\ &= \sum_a E^a_a
    \otimes q^{h_{\varepsilon_a}} \otimes q^{h_{\varepsilon_a}} \notag
    \\ & \qquad + (q-q^{-1}) \sum_{\varepsilon_{b} > \varepsilon_{a}}
    (-1)^{[b]} E^a_b \otimes \bigl( q^{h_{\varepsilon_a}}
    \hat{\sigma}_{ba} \otimes q^{h_{\varepsilon_a}} +
    q^{h_{\varepsilon_b}} \otimes q^{h_{\varepsilon_a}}
    \hat{\sigma}_{ba} \bigr) \notag \\ & \qquad + (q-q^{-1})^2
    \sum_{\varepsilon_{b} > \varepsilon_{c} >\varepsilon_a}
    (-1)^{[b]+[c]} E^a_b \otimes q^{h_{\varepsilon_c}}
    \hat{\sigma}_{bc} \otimes q^{h_{\varepsilon_a}}
    \hat{\sigma}_{ca}. \label{R13R12}
\end{align*}

\noindent Also, the coproduct properties \eqref{coprod} imply

\begin{multline*} 
(\text{id} \otimes \Delta) R = \sum_a E^a_a \otimes
  q^{h_{\varepsilon_a}} \otimes q^{h_{\varepsilon_a}} \\+ (q-q^{-1})
  \sum_{\varepsilon_{b} > \varepsilon_{a}} (-1)^{[b]} E^a_b \otimes
  (q^{h_{\varepsilon_a}} \otimes q^{h_{\varepsilon_a}})
  \Delta(\hat{\sigma}_{ba}).
\end{multline*}

\noindent Hence $R$ will satisfy equation \eqref{1tensordel} if and
only if $\Delta(\hat{\sigma}_{ba})$ is given by:

\begin{equation*}
\Delta(\hat{\sigma}_{ba}) = \hat{\sigma}_{ba} \otimes I +
q^{h_{\varepsilon_b} - h_{\varepsilon_a}} \otimes \hat{\sigma}_{ba} +
(q-q^{-1}) \!  \sum_{\varepsilon_{b} > \varepsilon_{c}
  >\varepsilon_a}\! (-1)^{[c]} q^{h_{\varepsilon_c} -
  h_{\varepsilon_a}} \hat{\sigma}_{bc} \otimes \hat{\sigma}_{ca}.
\end{equation*}

The simple values of $\hat{\sigma}_{ba}$ and the inductive relations
can be used to calculate $\Delta(\hat{\sigma}_{ba})$, and show that it
is indeed of this form.  First the coproduct of each of the simple
operators given in equation \eqref{ivalues} is directly calculated,
and found to be of the form required.  To find the coproduct for the
remaining values of $\hat{\sigma}_{ba}$ we use the inductive relations
\eqref{indrel}:

\begin{equation*}
\hat{\sigma}_{ba} = q^{-(\varepsilon_b,\varepsilon_a)}
  \hat{\sigma}_{bc}
  \hat{\sigma}_{ca}-q^{-(\varepsilon_c,\varepsilon_c)}(-1)^{([b]+[c])([a]+[c])}
  \hat{\sigma}_{ca} \hat{\sigma}_{bc}, \quad \varepsilon_b >
  \varepsilon_c > \varepsilon_a,
\end{equation*}

\noindent where $c \neq \overline{b}$ or $\overline{a}$.  We assume
our formula for the coproduct holds for $\hat{\sigma}_{bc}$ and
$\hat{\sigma}_{ca}$, where $\varepsilon_b > \varepsilon_c >
\varepsilon_a$, and then show it is also true for $\hat{\sigma}_{ba}$.

The coproduct is an algebra homomorphism, so for $\varepsilon_b >
\varepsilon_c > \varepsilon_a$, $c \neq \overline{a}$ or $
\overline{b}$, we have

\begin{equation*}
\Delta(\hat{\sigma}_{ba}) = q^{-(\varepsilon_b,\varepsilon_a)}
  \Delta(\hat {\sigma}_{bc}) \Delta(\hat{\sigma}_{ca})
  -q^{-(\varepsilon_c, \varepsilon_c)} (-1)^{([b]+[c])([a]+[c])}
  \Delta(\hat{\sigma}_{ca})\Delta(\hat{\sigma}_{bc}).
\end{equation*}

\noindent We can always choose $c$ satisfying the conditions such that either
$\varepsilon_b - \varepsilon_c$ or $\varepsilon_c - \varepsilon_a$ is
a simple root.  By considering each of those cases separately,
substituting in the expressions for $\Delta(\hat{\sigma}_{bc})$ and
$\Delta(\hat{\sigma}_{ca})$ and using the relations \eqref{qcom} and
\eqref{indrel} to simplify, after lengthy but straightforward
calculations given in \cite{me} we obtain

\begin{equation*}
\Delta(\hat{\sigma}_{ba}) = \hat{\sigma}_{ba} \otimes I + q^{h_{ba}}
  \otimes \hat{\sigma}_{ba} + (q-q^{-1}) \sum_{\varepsilon_{b} >
  \varepsilon_{c} > \varepsilon_a}(-1)^{[c]}
  q^{h_{ca}}\hat{\sigma}_{bc} \otimes\hat{\sigma}_{ca}
\end{equation*}

\noindent as required.  As all the relations in Tables \ref{list},
\ref{even} and \ref{odd} have either $\varepsilon_b - \varepsilon_c$
or $\varepsilon_c - \varepsilon_a$ as a simple root, this is
sufficient to prove the formula for $\Delta(\hat{\sigma}_{ba})$ for
all $\varepsilon_b > \varepsilon_a$.  As a check, however, we also
confirmed that this is consistent with the $q$-commutation relations
\eqref{qcom}, although again the calculations are tedious and omitted.
Hence we have proven that the matrix $R$ satisfies the property

\begin{equation*}
(\text{id} \otimes \Delta) R = R_{13} R_{12}.
\end{equation*}


\subsection{The Intertwining Property}

\noindent To confirm that we have a Lax operator we need to check one
last relation, namely the intertwining property for the other
generators.

\begin{equation} \label{intertwine}
R \Delta(a) = \Delta^T (a) R, \qquad \forall a \in U_q[osp(m|n)].
\end{equation}

\noindent Now $R$ is weightless, so it commutes with all the Cartan elements.
Moreover, $\Delta(q^{h_a}) = \Delta^T (q^{h_a}), \forall h_a \in H$,
so the Cartan elements will automatically satisfy equation
\eqref{intertwine}.  Thus it remains only to verify the intertwining
property for the lowering generators, $f_a$.  Unfortunately, knowing
the raising generators satisfy the intertwining property does not
appear helpful.  Instead, we start by assuming the form of the Lax
operator and that it satisfies the intertwining property for the
lowering generators, and then proceed as in Section \ref{Developing
Relations}.  Again two different equations are obtained, this time

\begin{multline}
f_c \otimes (q^{\frac{3}{2}h_c} - q^{-\frac{1}{2} h_c}) \\ = (q -
q^{-1}) \sum_{\varepsilon_{b} - \varepsilon_{a} = \alpha_c} (-1)^{[b]}
E^a_b \otimes \bigl( q^{\frac{1}{2} (\alpha_c, \varepsilon_b)}
\hat{\sigma}_{ba} f_c - (-1)^{[c]} q^{\frac{1}{2} (\alpha_c,
\varepsilon_a)} f_c \hat{\sigma}_{ba} \bigr) \label{basef}
\end{multline}

\noindent and

\begin{multline}
q^{\frac{3}{2} (\varepsilon_b -\varepsilon_a -\alpha_c, \alpha_c)}
  \langle a|f_c|a' \rangle \hat{\sigma}_{ba'} q^{\frac{3}{2} h_c} -
  (-1)^{([a]+[b])[c]} \langle b'|f_c|b \rangle \hat{\sigma}_{b'a}
  q^{-\frac{1}{2} h_c} \\ = q^{\frac{1}{2} (\alpha_c, \varepsilon_b)}
  \hat{\sigma}_{ba} f_c - (-1)^{([a]+[b])[c]} q^{\frac{1}{2}
  (\alpha_c, \varepsilon_a)} f_c \hat{\sigma}_{ba}, \quad
  \varepsilon_b > \varepsilon_a. \label{eqf}
\end{multline}

\

By considering the simple roots individually and using the
$U_q[osp(m|n)]$ defining relations, equation \eqref{basef} can be
shown to be consistent with the simple values of $\hat{\sigma}_{ba}$
given in Table \ref{fundval} \cite{me}, as expected.
 The various relations that can be deduced from equation \eqref{eqf}
 are no longer inductive in form, so can not be used to directly
 construct the $\hat{\sigma}_{ba}$ for direct comparison.  They can,
 however, be shown to be consistent with relations \eqref{qcom} and
 \eqref{indrel}, as demonstrated in \cite{me}.  Thus the conditions
 obtained by considering the intertwining property for the lowering
 generators,

\begin{equation*}
R\Delta(f_c) = \Delta^T(f_c)R
\end{equation*}

\noindent are satisfied by the matrix $R$, and hence $R$ satisfies the
intertwining property for all elements $a \in U_q[osp(m|n)]$.


\subsection{The Lax Operator}

\noindent We have now proven, as expected, that the matrix $R$ found
earlier satisfies both the intertwining property and $(\text{id}\,
\otimes\, \Delta)R = R_{13} R_{12}.$ The other $R$-matrix property,
containing $(\Delta\, \otimes\, \text{id})R$, is clearly not
applicable here.  It is not necessary, however, as we know there is a
Lax operator belonging to $\pi \bigl( U_q[osp(m|n)]^-\bigr) \otimes
U_q[osp(m|n)]^+$, and we have shown there is only one such
possibility.  Thus the work in the previous two sections confirms the
following proposition:

\begin{proposition}

The Lax operator, $R = (\pi \otimes \text{id}) \mathcal{R}$ for the
quantum superalgebra $U_q[osp(m|n)]$, where $\mathcal{R} \in
U_q[osp(m|n)]^- \otimes U_q[osp(m|n)]^+$ and $m >2$, is given by

\begin{align*}
R &= q ^ {h_{x} \otimes h^{x}} \Bigl[ I \otimes I + (q - q^{-1})
  \sum_{\varepsilon_{a} < \varepsilon_{b}} (-1)^{[b]} E^a_b \otimes
  \hat{\sigma}_{ba} \Bigr] \\ &= \sum_a E^a_a \otimes
  q^{h_{\varepsilon_a}} + (q-q^{-1}) \sum_{\varepsilon_a <
  \varepsilon_b} (-1)^{[b]} E^a_b \otimes q^{h_{\varepsilon_a}}
  \hat{\sigma}_{ba},
\end{align*}

\noindent where the operators $\hat{\sigma}_{ba}$ satisfy equations
\eqref{ivalues}, \eqref{qcom} and \eqref{indrel}.

\end{proposition}

As an aside, the two properties verified directly are sufficient to
prove $R$ satisfies the Yang-Baxter equation.  For using only those,
we see

\begin{align*}
R_{23} R_{13} R_{12} &= R_{23} (\text{id} \otimes \Delta)R \\ &=
[(\text{id} \otimes \Delta ^T )R] R_{23} \\ &= [(\text{id} \otimes T)
((\text{id} \otimes \Delta)R] R_{23} \\ &= [(\text{id} \otimes T)
R_{13} R_{12}] R_{23} \\ &= R_{12} R_{13} R_{23}
\end{align*}

\noindent as required.  It is very surprising that there is a unique solution to

\begin{equation*}
R \Delta (e_c) = \Delta^T (e_c) R,
\end{equation*}

\noindent given we only considered elements of $\pi
\bigl( U_q[osp(m|n)]^-\bigr) \otimes U_q[osp(m|n)]^+$. 

\subsection{The Opposite Lax Operator} \label{opposite}

\noindent Having found the Lax operator $R = (\pi \otimes \text{id})
\mathcal{R}$, we wish to use that result to find its opposite $R^T =
(\pi \otimes \text{id}) \mathcal{R}^T$, where $\mathcal{R}^T$ is the
opposite universal $R$-matrix of $U_q[osp(m|n)]$.  We begin by showing
that $\mathcal{R}^T$ is in fact equal to $\mathcal{R}^\dagger$, where
$^\dagger$ represents graded conjugation, defined below.

A graded conjugation on $U_q[osp(m|n)]$ is defined on the simple
generators by:

\begin{equation*}
e_a^\dagger = f_a, \qquad f_a^\dagger = (-1)^{[a]} e_a, \qquad
h_a^\dagger = h_a.
\end{equation*}

\noindent It is consistent with the coproduct and extends naturally to
all remaining elements of $U_q[osp(m|n)]$, satisfying the properties:

\begin{align*}
&(\sigma^a_b)^\dagger = (-1)^{[a]([a]+[b])} \sigma^b_a,\\
&(ab)^\dagger = (-1)^{[a][b]} b^\dagger a^\dagger, \\ &(a \otimes
b)^\dagger = a^\dagger \otimes b^\dagger, \\ &\Delta (a)^\dagger =
\Delta (a^\dagger).
\end{align*}

Returning to the universal $R$-matrix $\mathcal{R}$, we know

\begin{alignat*}{3}
&&&\mathcal{R} \Delta(a) = \Delta^T(a) \mathcal{R}, &\qquad &\forall a
\in U_q [osp(m|n)], \notag \\ &\Rightarrow& \quad & \Delta(a)^\dagger
\mathcal{R}^\dagger = \mathcal{R}^ \dagger\Delta^T(a)^\dagger \notag
\\ &\Rightarrow && \Delta(a^\dagger) \mathcal{R}^\dagger =
\mathcal{R}^\dagger \Delta^T(a^\dagger) \notag \\ &\Rightarrow &&
\Delta(a) \mathcal{R}^\dagger = \mathcal{R}^\dagger \Delta^T(a),&&
\forall a \in U_q[osp(m|n)].
\end{alignat*}

Similarly, $\mathcal{R}^\dagger$ satisfies the other
$R$-matrix properties \eqref{Requations}.  As there is a unique
universal $R$-matrix belonging to $U_q[osp(m|n)]^+ \otimes
U_q[osp(m|n)]^-$, the only possibility is $\mathcal{R}^T =
\mathcal{R}^\dagger$.

Now it is known that the vector representation is
superunitary.  A discussion of superunitary representations is given
in \cite{LG}, where they are called grade star representations, but
here we need only note this implies

\begin{equation*}
\pi (a^\dagger) = \pi (a)^\dagger, \qquad \forall a \in U_q[osp(m|n)].
\end{equation*}

\noindent Hence

\begin{align*}
R^T &= (\pi \otimes \text{id}) \mathcal{R}^\dagger \notag \\ &= [(\pi
\otimes \text{id}) \mathcal{R}]^\dagger \notag \\ &= R^\dagger.
\end{align*}

Thus we can find the opposite Lax operator $R^T$ simply by
using the usual rules for graded conjugation.  As $R$ is given by

\begin{equation*}
R = \sum_a E^a_a \otimes q^{h_{\varepsilon_a}} + (q-q^{-1})
\sum_{\varepsilon_b > \varepsilon_a} (-1)^{[b]} E^a_b \otimes
q^{h_{\varepsilon_a}} \hat{\sigma}_{ba},
\end{equation*}

\noindent we find the opposite Lax operator $R^T$ can be written as

\begin{equation} \label{RT}
R^T = \sum_a E^a_a \otimes q^{h_{\varepsilon_a}} + (q-q^{-1})
\sum_{\varepsilon_b > \varepsilon_a} (-1)^{[a]} E^b_a \otimes
\hat{\sigma}_{ab} q^{h_{\varepsilon_a}},
\end{equation}

\noindent where

\begin{equation*}
\hat{\sigma}_{ab} = (-1)^{[b]([a]+[b])} \hat{\sigma}_{ba}^\dagger,
\quad \varepsilon_b > \varepsilon_a.
\end{equation*}


\subsection{q-Serre Relations}
\noindent Having shown the relations found in Section \ref{noniv}
define a Lax operator, we also wish to see if they incorporate the
$q$-Serre relations.  It is too space-consuming to list all of these,
so we will merely provide a couple of examples, including the extra
$q$-Serre relations.

First recall that if $\varepsilon_b - \varepsilon_a$ is a
simple root, then $\hat{\sigma}_{ba} \propto e_c q^{\frac{1}{2}h_c}$
for either $c = b$ or $c=\overline{a}$.  Then setting $E_a = e_a
q^{\frac{1}{2} h_a}$, we see from Definition \ref{def} that:

\begin{alignat}{2}
&&\Delta (E_a) &= q^{h_a} \otimes E_a + E_a \otimes 1 \notag \\
&&S(E_a) &= -q^{-\frac{1}{2}(\alpha_a, \alpha_a)}
q^{-\frac{1}{2}h_a}e_a \notag \\ &&&= - q^{-h_a} E_a \notag \\
&\therefore \qquad & ad\, E_a \circ b &= - (-1)^{[a][b]} q^{h_a}b
q^{-h_a}E_a + E_a b \notag \\ &&&= E_a b - (-1)^{[a][b]}
q^{(\alpha_a,\varepsilon_b)} bE_a. \label{adjoint}
\end{alignat}

\noindent Now consider the simple generators $\hat{\sigma}_{i\,i+1}$
and $\hat{\sigma}_{i+1\, i+2}$.

\begin{align*}
(ad\: \hat{\sigma}_{i\,i+1}\: \circ)^2 \hat{\sigma}_{i+1\, i+2} &=
  ad\: \hat {\sigma}_{i\,i+1} \circ (\hat{\sigma}_{i\,i+1}
  \hat{\sigma}_{i+1\, i+2} - q^{-1} \hat{\sigma}_{i+1\, i+2}
  \hat{\sigma}_{i\,i+1}) \\ &= ad\: \hat{\sigma}_{i\,i+1} \circ
  \hat{\sigma}_{i\,i+2} \\ &= \hat{\sigma}_{i\,i+1}
  \hat{\sigma}_{i\,i+2} - q \hat{\sigma}_{i\,i+2}
  \hat{\sigma}_{i\,i+1} \\ &= 0
\end{align*}

\noindent from \eqref{qcom}.  This is equivalent to the $q$-Serre
relation $(ad\: e_b\: \circ)^{1-a_{bc}} e_c = 0$ for this pair of
simple operators.  In a similar way, we can verify this relation for
any $b \neq c, \, (\varepsilon_b, \varepsilon_b) \neq 0$.  
The defining relations for the $\hat{\sigma}_{ba}$,
therefore, incorporate all the standard $q$-Serre relations for
raising generators.

This still leaves the extra $q$-Serre relations, which
involve the odd root.  There are only two of these for our choice of
simple roots \cite{Yamane}.  Explicitly, taking into account the
different conventions between \cite{Yamane} and here, the relevant extra 
$q$-Serre relations for $U_q[osp(m|n)]$ can be written as

\begin{align} \label{q1}
&\bigl[ \hat{\sigma}_{\mu=k\,i=1}, \bigl[
\hat{\sigma}_{\nu=k-1\,\mu=k},
[\hat{\sigma}_{\mu=k\,i=1},\hat{\sigma}_{i=1\,j=2} ]_q \, ]_q \,] = 0
\\ &\bigl[ \hat{\sigma}_{\mu=k\,i=1}, \bigl[ \hat{\sigma}_{i=1\,j=2},
[\hat{\sigma}_{\mu=k\,i=1}, \hat{\sigma}_{\nu=k-1\,\mu=k}]_q \, ]_q
\,] = 0, \label{q2}
 \end{align}

\noindent where $[x,y]_q$ represents the adjoint action $ad\, x \circ
y$.

Consider equation \eqref{q1}.  Using the defining relations
\eqref{qcom} and \eqref{indrel} for the $\hat{\sigma}_{ba}$ together
with the adjoint action as given in equation \eqref{adjoint}, we find:

\begin{align*}
&[\hat{\sigma}_{\mu=k\,i=1}, [ \hat{\sigma}_{\nu=k-1\,\mu=k},
  [\hat{\sigma}_{\mu=k\,i=1},\hat{\sigma}_{i=1\,j=2} ]_q \, ]_q
  \,]\notag\\ =\:&[ \hat{\sigma}_{\mu=k\,i=1}, [
  \hat{\sigma}_{\nu=k-1\,\mu=k}, (\hat{\sigma}_{\mu=k\,i=1}
  \hat{\sigma}_{i=1\,j=2} - q^{-1} \hat{\sigma}_{i=1\,j=2}
  \hat{\sigma}_{\mu=k\,i=1}) ]_q \,]\notag\\ =\:&[
  \hat{\sigma}_{\mu=k\,i=1}, [ \hat{\sigma}_{\nu=k-1\,\mu=k},
  \hat{\sigma}_{\mu=k\, j=2} ]_q \,] \notag \\ =\:&[
  \hat{\sigma}_{\mu=k\,i=1}, (\hat{\sigma}_{\nu=k-1\,\mu=k}
  \hat{\sigma}_{\mu=k\, j=2} - q \hat{\sigma}_{\mu=k\, j=2}
  \hat{\sigma}_{\nu=k-1\,\mu=k})] \notag \\ =\:&[
  \hat{\sigma}_{\mu=k\,i=1}, \hat{\sigma}_{\nu=k-1\,j=2}] \notag \\
  =\:& 0
\end{align*}

\noindent as required.  It is equally straightforward to show that
equation \eqref{q2} arises from the defining relations of the
$\hat{\sigma}_{ba}$.  Hence these compact defining relations for the
$\hat{\sigma}_{ba}$ incorporate not only the standard $q$-Serre
relations for the raising generators, but also the extra ones.  This
is quite remarkable, as the equivalent $q$-Serre relations are not
used in the derivation.


\section{The $R$-matrix for the Vector Representation} \label{vector}

\noindent The Lax operator can be used to explicitly calculate an
$R$-matrix for any tensor product representation $\pi \otimes \phi$,
where $\phi$ is an arbitrary representation.  In particular, it
provides a more straightforward method of calculating $R$ for the
tensor product of the vector representation, $\pi \otimes \pi$, than
using the tensor product graph method \cite{Mehta}.

By specifically constructing the $R$-matrix for the vector
representation, we also illustrate concretely the way the recursion
relations can be applied to find the $R$-matrix for an arbitrary
representation.  Although the values for $\hat{\sigma}_{ba}$ obtained
will change for each representation, they can always be constructed by
applying the same equations in the same order.  We could choose to use
only the relations listed in the tables in the appendix, but using the
general form of the inductive relations shortens and simplifies the
process.  Only the method of calculation is included here, but the
full calculations can be found in \cite{me}.

First the vector representation is applied to the simple
operators given in Table \ref{fundval}, with the results written below
in Table \ref{basevector}.

\begin{table}[ht] 
\caption{The simple values of $\hat{\sigma}_{ba}$ in the vector
representation.}\label{basevector}
\centering
\begin{tabular}{|l|lll|} 
\hline\noalign{\smallskip}
\multicolumn{1}{|c|}{Simple Root}& \multicolumn{3}{c|}{Corresponding
  $\hat{\sigma}_{ba}$} \\ 
\noalign{\smallskip}\hline\noalign{\smallskip}
$\alpha_i = \varepsilon_i -
  \varepsilon_{i+1},\,i < l$ & $\hat{\sigma}_{i\, i+1}$&$ =
  -\hat{\sigma}_{\overline{i+1}\,\overline{i}}$&$ = E^i_{i+1} - E^
  {\overline{i+1}}_{\overline{i}}$ \\ $\alpha_l = \varepsilon_{l-1} +
  \varepsilon_l,\,m = 2l$ & $\hat{\sigma}_ {l-1\,\overline{l}}$&$ =
  -\hat{\sigma}_{l\, \overline{l-1}}$&$ = E^{l-1}_ {\overline{l}} -
  E^l_{\overline{l-1}}$ \\ $\alpha_l = \varepsilon_l,\,m=2l+1$ &
  $\hat{\sigma}_{l\, l+1}$&$ = -q^{-\frac {1}{2}} \hat{\sigma}_{l+1\,
  \overline{l}}$&$ = E^l_{l+1} - q^{-\frac{1}{2}}
  E^{l+1}_{\overline{l}}$ \\ $\alpha_\mu = \delta_\mu -
  \delta_{\mu+1},\,\mu < k$ & $\hat{\sigma}_ {\mu\, \mu+1} $&$ =
  \hat{\sigma}_{\overline{\mu+1}\, \overline{\mu}} $&$= E^\mu_{\mu+1}
  + E^{\overline{\mu+1}}_{\overline{\mu}}$ \\ $\alpha_s = \delta_k -
  \varepsilon_1,$ & $\hat{\sigma}_{\mu=k\, i=1} $&$ = (-1)^k q
  \hat{\sigma}_{\overline{i}=\overline{1}\, \overline{\mu} =
  \overline{k}} $&$= E^{\mu=k}_{i=1} + (-1)^k q E^{\overline{i} =
  \overline{1}}_{\overline{\mu} = \overline{k}}$ \\ 
\noalign{\smallskip}\hline
\end{tabular}
\end{table}

Then the inductive relations \eqref{indrel} are applied to
these simple operators to find the remaining values of
$\hat{\sigma}_{ba}$.  One of many equivalent ways of doing this is
given below.  

Construct
\begin{enumerate}
\item $\hat{\sigma}_{ji},\hat{\sigma}_{\overline{i}\, \overline{j}}$
  for $1 \leq j < i \leq \lceil \frac{m}{2} \rceil$, using
  $\hat{\sigma}_{i\, i+1},
  \hat{\sigma}_{\overline{i+1}\,\overline{i}}$ for $1 \leq i <l$ and
  $\hat{\sigma}_{l\, l+1},\hat{\sigma}_{l+1\, \overline{l}}$ when
  $m=2l+1$
  \label{item ji}
\item $\hat{\sigma}_{\nu\mu},\hat{\sigma}_{\overline{\mu}\,
  \overline{\nu}}$ for $1 \leq \nu < \mu \leq k$, using
  $\hat{\sigma}_{\mu\, \mu+1},\hat{\sigma} _{\overline{\mu+1}\,
  \overline{\mu}}$ for $1 \leq \mu < k$ \label{item numu}
\item $\hat{\sigma}_{\mu i}, \hat{\sigma}_{\overline{i}\,
  \overline{\mu}}$ for $1 \leq i \leq \lceil \frac{m}{2} \rceil, 1
  \leq \mu \leq k$, using $\hat{\sigma}_{\mu=k\, i=1}$ and
  $\hat{\sigma}_{\overline{i}=\overline{1}\, \overline{\mu}
  =\overline{k}} $ together with the values calculated in steps
  \ref{item ji} and \ref{item numu} \label{item mui}
\item $\hat{\sigma}_{i\, \overline{j}}$ for $1 \leq i,j \leq l$, using
  $\hat{\sigma}_{l\, l+1}$ and $\hat{\sigma}_{l+1\, \overline{l}}$
  when $m=2l+1$ or $\hat{\sigma}_{l-1\, \overline{l}},
  \hat{\sigma}_{l\, \overline {l-1}}$ and $\hat{\sigma}_{l
  \overline{l}}=0$ when $m=2l$, together with the results from step
  \ref{item ji} \label{item ijbar}
\item $\hat{\sigma}_{i\,\overline{\mu}}, \hat{\sigma}_{\mu\,
  \overline{i}}$ for $1 \leq i \leq l, 1 \leq \mu \leq k$, using the
  results from steps \ref{item mui} and \ref{item ijbar} \label{item
  imubar}
\item $\hat{\sigma}_{\mu\,\overline{\nu}}$ for $1 \leq \mu,\nu \leq
  k$, using the results from steps \ref{item mui} and \ref{item
  imubar}
\end{enumerate}

Following this procedure in the vector representation, we
find the following form for the operators $\hat{\sigma}_{ba},\;
\varepsilon_b > \varepsilon_a$:

\begin{equation*}
\hat{\sigma}_{ba} = q^{-(\varepsilon_a, \varepsilon_b)} E^b_a -
(-1)^{[b]([a]+[b])} \xi_a \xi_b q^{(\varepsilon_a, \varepsilon_a)}
q^{(\rho, \varepsilon_a - \varepsilon_b)}
E^{\overline{a}}_{\overline{b}}, \qquad \varepsilon_b > \varepsilon_a.
\end{equation*}

Thus we have shown the $R$-matrix for the vector
rep. of $U_q[osp(m|n)]$, $\mathfrak{R} = (\pi \otimes \pi)
\mathcal{R}$, is given by

\begin{equation*}
\mathfrak{R} = q ^{h_{j} \otimes h^{j}} \Bigl[ I \otimes I + (q - q^{-1})
\sum_{\varepsilon_{b} > \varepsilon_{a}} (-1)^{[b]} E^a_b \otimes
\hat{\sigma}_{ba} \Bigr]
\end{equation*}

\noindent where

\begin{equation*}
\hat{\sigma}_{ba} = q^{-(\varepsilon_a, \varepsilon_b)} E^b_a -
(-1)^{[b]([a]+[b])} \xi_a \xi_b q^{(\varepsilon_a, \varepsilon_a)}
q^{(\rho, \varepsilon_a - \varepsilon_b)}
E^{\overline{a}}_{\overline{b}}.
\end{equation*}

\noindent This can be written in a more elegant form, namely

\begin{equation*}
\mathfrak{R} = \sum_{a,b} q^{(\varepsilon_{a}, \varepsilon_{b})} E^a_a \otimes
E^b_b \; + \; (q-q^{-1}) \sum_{\varepsilon_{b} > \varepsilon_{a}}
(-1)^{[b]} E^a_b \otimes \tilde{\sigma}_{ba}
\end{equation*}

\noindent where $\tilde{\sigma}_{ba} = q^{h_{\varepsilon_a}}
\hat{\sigma}_{ba}$.  Hence we have the following result:

\begin{proposition}  The $R$-matrix for the vector representation, $\mathfrak{R} = (\pi \otimes \pi) \mathcal{R}$, is given by 

\begin{equation*}
\mathfrak{R} = \sum_{a,b} q^{(\varepsilon_{a}, \varepsilon_{b})} E^a_a \otimes
E^b_b \; + \; (q-q^{-1}) \sum_{\varepsilon_{b} > \varepsilon_{a}}
(-1)^{[b]} E^a_b \otimes \tilde{\sigma}_{ba},
\end{equation*} 

\noindent where

\begin{equation*} 
\tilde{\sigma}_{ba} = E^b_a - (-1)^{[b]([a]+[b])} \xi_a \xi_b
q^{(\rho, \varepsilon_a - \varepsilon_b)}
E^{\overline{a}}_{\overline{b}}, \qquad \varepsilon_b > \varepsilon_a.
\end{equation*}
\end{proposition}

\

 We can also explicitly find the opposite $R$-matrix $\mathfrak{R}^T$ as
given in equation \eqref{RT}, using

\begin{equation*}
(E^a_b)^\dagger = (-1)^{[a]([a]+[b])} E^b_a.
\end{equation*}

\noindent We obtain this result:

\begin{proposition}
\noindent The opposite $R$-matrix for the vector representation, $\mathfrak{R}^T
= (\pi \otimes \pi) \mathcal{R}^T$, is given by

\begin{equation*}
\mathfrak{R}^T = \sum_{a,b} q^{(\varepsilon_a, \varepsilon_b)} E^a_a \otimes
E^b_b + (q-q^{-1}) \sum_{\varepsilon_b > \varepsilon_a} (-1)^{[a]}
E^b_a \otimes \tilde{\sigma}_{ab},
\end{equation*}

\noindent where

\begin{equation*}
\tilde{\sigma}_{ab} = E^a_b - (-1)^{[a]([a]+[b])} \xi_a \xi_b
  q^{(\rho, \varepsilon_a - \varepsilon_b)}
  E^{\overline{b}}_{\overline{a}}.
\end{equation*}

\end{proposition}

\

These formulae for $\mathfrak{R}$ and $\mathfrak{R}^T$ on the vector
representation agree with those given in \cite{Mehta}.  In that thesis
the $R$-matrix for the vector representation was calculated using
projection operators onto invariant submodules of the tensor product.
The greatest advantage of the current method is it gives a
straightforward way of constructing a solution to the Yang-Baxter
Equation in an arbitrary representation of $U_q[osp(m|n)]$.  

Applying the tensor product graph method to $\mathfrak{R}$, it was
shown in \cite{Mehta} that the spectral dependent $R$-matrix for the
vector representation of the quantum affine superalgebra
$U_q[osp(m|n)^{(1)}]$ is

{\small
\begin{multline*}
\mathfrak{R}(z)  =  \frac{(q-q^{-1})z} {(q-zq^{-1})}P - \frac{(q-q^{-1})z(z-1)}
  {(q-q^{-1}z)(z-q^{m-n-2})} \sum_{a,b} (-1)^{[a][b]} \xi_a \xi_b q^{(\rho, 
  \varepsilon_a-\varepsilon_b)} E^a_b \otimes E^{\overline{a}}_{\overline{b}}\\
-\frac{(z-1)}{(q -zq^{-1})} \biggl\{ I + (q^{1/2}-q^{-1/2}) \sum_a (-1)^{[a]} 
  E^a_a \otimes \hat{\sigma}^a_a + (q-q^{-1})\sum_{\varepsilon_a<\varepsilon_b}
  (-1)^{[b]} E^a_b \otimes \hat{\sigma}^b_a \biggr\}
\end{multline*}}

\noindent and the spectral dependent $R$-matrix for the vector representation 
 of the quantum twisted affine superalgebra $U_q[gl(m|n)^{(2)}]$ is

{\small
\begin{multline*}
\mathfrak{R}(z) = \frac{(q-q^{-1}) z}{(q-zq^{-1})} P - \frac{(q-q^{-1})z(z-1)} 
  {(q-q^{-1}z)(z+q^{m-n})} \sum_{a,b} (-1)^{[a][b]} \xi_a \xi_b
  q^{(\rho, \varepsilon_a-\varepsilon_b)} E^a_b \otimes
  E^{\overline{a}}_{\overline{b}} \\ 
-\frac{(z-1)}{(q -zq^{-1})} \biggl\{ I + (q^{1/2}-q^{-1/2}) \sum_a (-1)^{[a]} 
  E^a_a \otimes \hat{\sigma}^a_a +(q-q^{-1}) \sum_{\varepsilon_a<\varepsilon_b}
  (-1)^{[b]} E^a_b \otimes \hat{\sigma}^b_a \biggr\}.
\end{multline*}}

In a similar way the tensor product graph can be applied to the Lax Operator 
to find the spectral dependent $R$-matrix for other affinisable representations
 of the form $\pi \otimes \sigma$ applied to $U_q[osp(m|n)^{(1)}]$ or 
$U_q[gl(m|n)^{(2)}]$.


\section{Conclusion} 

\noindent In this Communication a Lax operator was constructed for the
 $B$ and $D$ series of quantum superalgebras. A general ansatz for the
 Lax operator in terms of unknown elements of $U_q[osp(m|n)]$ was
 assumed.  Formulae identifying the simple operators were found, and
 then a set of inductive and $q-$commutative relations developed that
 can be used to calculate the remaining non-simple operators. This
 result is universal and thus the Lax operator generates a solution of
 the Yang-Baxter equation on the space $V\otimes V\otimes W$ for any
 module $W$.  A specific example was given in Section \ref{vector},
 where the $R$-matrix for the vector representation was calculated
 from the Lax operator.  Together with the results of \cite{Zhang2}
 for the $A$ series and \cite{CLax} for the $C$ series, this completes
 the construction of Lax operators for all non-exceptional quantum
 superalgebras.

\vspace{0.5cm}
\noindent{\bf Acknowledgements --} We gratefully acknowledge financial
support from the Australian Research Council.

\newpage
\appendix

\section{The Relations Governing the Operators $\hat{\sigma}_{ba}$}

\begin{table}[ht]
\caption{The relations for the operators $\hat{\sigma}_{ba}$ common to
all values of $m$} 
\label{list} 
\centering
\begin{tabular}{|ll|} 
\hline\noalign{\smallskip}
\hspace{25mm} Relation & \hspace{10mm} Conditions \\ 
\noalign{\smallskip}\hline\noalign{\smallskip}
$\hat{\sigma}_{b\, i+1} = \hat{\sigma}_{bi} \hat{\sigma}_{i\, i+1} -
  q^{-1} \hat{\sigma}_{i\, i+1} \hat{\sigma}_{bi},$&$i<l,\:
  \varepsilon_b > \varepsilon_i$ \\ $\hat{\sigma}_{\overline{i+1}\,a}
  = \hat{\sigma}_{\overline{i+1}\,\overline{i}}
  \hat{\sigma}_{\overline{i}\, a} - q^{-1}
  \hat{\sigma}_{\overline{i}\, a}
  \hat{\sigma}_{\overline{i+1}\,\overline{i}}, $&$ i<l,\:
  \varepsilon_a < - \varepsilon_i$\\ $\hat{\sigma}_{b\,\overline{i}} =
  q^{(\alpha_i, \varepsilon_b)} \hat{\sigma}_ {b\, \overline{i+1}}
  \hat{\sigma}_{\overline{i+1}\, \overline{i}} - q^{-1}
  \hat{\sigma}_{\overline{i+1}\,\overline{i}}
  \hat{\sigma}_{b\,\overline{i+1}}, $&$ i<l,\: \varepsilon_b > -
  \varepsilon_{i+1},$ \\ &$ b \neq i+1$ \\ $\hat{\sigma}_{ia} =
  q^{-(\alpha_i, \varepsilon_a)} \hat{\sigma}_{i\, i+1}
  \hat{\sigma}_{i+1\, a} - q^{-1} \hat{\sigma}_{i+1\, a}
  \hat{\sigma}_{i\,i+1}, $&$ i<l,\: \varepsilon_a <
  \varepsilon_{i+1},$ \\&$ a \neq \overline{i+1}$ \\
  $\hat{\sigma}_{i+1\,\overline{i}} + \hat{\sigma}_{i\,
  \overline{i+1}} = q^{-1} [\hat{\sigma}_{i\, i+1},
  \hat{\sigma}_{i+1\, \overline{i+1}}],$& $i<l$\\
  $q^{(\alpha_i,\varepsilon_b)} \hat{\sigma}_{ba} \hat{\sigma}_{i\,
  i+1} - q^{-(\alpha_i,\varepsilon_a)} \hat{\sigma}_{i\, i+1}
  \hat{\sigma}_{ba} = 0,$ & $i<l;\: \varepsilon_b >
  \varepsilon_a;$\\&$ a \neq i, \overline{i+1}$ \\ & and $b \neq i+1,
  \overline{i}$ \\
$\hat{\sigma}_{\nu\, \mu+1} = \hat{\sigma}_{\nu \mu}
\hat{\sigma}_{\mu\, \mu+1} - q \hat{\sigma}_{\mu\, \mu+1}
\hat{\sigma}_{\nu \mu}, $&$ \nu < \mu<k$ \\
$\hat{\sigma}_{\overline{\mu+1}\, \overline{\nu}} =
\hat{\sigma}_{\overline {\mu+1} \, \overline{\mu}}
\hat{\sigma}_{\overline{\mu} \overline{\nu}} - q
\hat{\sigma}_{\overline {\mu} \overline{\nu}}
\hat{\sigma}_{\overline{\mu+1} \,\overline{\mu}}, $&$ \nu < \mu < k$
\\ $\hat{\sigma}_{b \overline{\mu}} = q^{(\alpha_\mu, \varepsilon_b)}
\hat{\sigma} _{b\, \overline{\mu+1}} \hat{\sigma}_{\overline{\mu+1}\,
  \overline{\mu}} - q \hat{\sigma}_{\overline{\mu+1}\, \overline{\mu}}
\hat{\sigma}_{b\, \overline {\mu+1}}, $&$ \mu<k,\: \varepsilon_b >
-\delta_{\mu + 1},$\\&$b\neq \mu +1$\\ $\hat{\sigma}_{\mu a} =
q^{-(\alpha_\mu, \varepsilon_a)} \hat{\sigma}_{\mu\, \mu+1}
\hat{\sigma}_{\mu + 1\, a} - q \hat{\sigma}_{\mu + 1\, a}
\hat{\sigma}_{\mu\, \mu+1},$&$ \mu<k,\: \varepsilon_a < \delta_{\mu +
  1},$\\ &$a \neq \overline{\mu + 1}$ \\ $\hat{\sigma}_{\mu+1\,
  \overline{\mu}} - \hat{\sigma}_{\mu\,\overline{\mu +1}} = q [
  \hat{\sigma}_{\mu +1 \, \overline{\mu
      +1}},\hat{\sigma}_{\mu\,\mu+1}],$ & $\mu<k$ \\
$q^{(\alpha_\mu,\varepsilon_b)} \hat{\sigma}_{ba} \hat{\sigma}_{\mu \,
  \mu +1} - q^{-(\alpha_\mu,\varepsilon_a)} \hat{\sigma}_{\mu \, \mu
  +1} \hat{\sigma}_ {ba} = 0, $&$\mu<k;\: \varepsilon_b >
\varepsilon_a;$ \\ &$a\neq \mu, \overline{\mu+1}$ \\ &and $b \neq
\mu+1, \overline{\mu}$ \\
$\hat{\sigma}_{\nu\, i=1} = \hat{\sigma}_{\nu\, \mu=k}
  \hat{\sigma}_{\mu=k \, i=1} - q \hat{\sigma}_{\mu=k \, i=1}
  \hat{\sigma}_{\nu\, \mu=k},$&$\nu <k$\\
  $\hat{\sigma}_{i=\overline{1}\, \overline{\nu}} =
  \hat{\sigma}_{i=\overline{1} \, \mu = \overline{k}}
  \hat{\sigma}_{\mu = \overline{k}\, \overline{\nu}} - q
  \hat{\sigma}_{\mu = \overline{k}\, \overline{\nu}} \hat{\sigma}_{i =
  \overline{1}\, \mu = \overline{k}}, $&$ \nu < k $\\
  $\hat{\sigma}_{\mu=k\, a} = q^{-(\alpha_s, \varepsilon_a)}
  \hat{\sigma}_{\mu=k \,i=1} \hat{\sigma}_{i=1\, a} - (-1)^{[a]}
  q^{-1} \hat{\sigma}_{i=1\, a} \hat {\sigma}_{\mu=k \,i=1},$&$
  \varepsilon_a < \varepsilon_1,\: a\neq i=\overline {1} $\\
  $\hat{\sigma}_{b\, \mu = \overline{k}} = q^{(\alpha_s,
  \varepsilon_b)} \hat{\sigma}_{b\, i= \overline{1}} \hat{\sigma}_{i =
  \overline{1}\, \mu = \overline{k}} - (-1)^{[b]} q^{-1}
  \hat{\sigma}_{i = \overline{1}\, \mu = \overline{k}}
  \hat{\sigma}_{b\, i= \overline{1}}, $&$ \varepsilon_b > -
  \varepsilon_1,\: b \neq i = 1$ \\ $\hat{\sigma}_{\mu=k\,
  \overline{i} = \overline{1}} - (-1)^k q \hat{\sigma}_ {i=1\,
  \overline{\mu} = \overline{k}} = q^{-1} [ \hat{\sigma}_{\mu=k\,
  i=1}, \hat{\sigma}_{i=1\, \overline{i} = \overline{1}}],$& \\
  $q^{(\alpha_s, \varepsilon_b)} \hat{\sigma}_{ba} \hat{\sigma}_{\mu=k
  \, i=1} - (-1)^{[a]+[b]} q^{-(\alpha_s,\varepsilon_a)}
  \hat{\sigma}_{\mu=k \, i=1} \hat{\sigma}_{ba} = 0 $&$ \varepsilon_b
  > \varepsilon_a;\: \varepsilon_a \neq \delta_k, -\varepsilon_1$ \\
  &and $\varepsilon_b \neq \varepsilon_1, -\delta_k$ \\ 
\noalign{\smallskip}\hline
\end{tabular}
\end{table}
\begin{table}[h]
\caption{The relations for the operators $\hat{\sigma}_{ba}$ that hold
only for even $m$} 
\label{even} 
\centering
\begin{tabular}{|ll|}
\hline\noalign{\smallskip}
\hspace{20mm} Relation & Conditions \\ 
\noalign{\smallskip}\hline\noalign{\smallskip}
$\hat{\sigma}_{b\, \overline{l-1}} = q^{(\alpha_l,\varepsilon_b)}
  \hat{\sigma} _{bl} \hat{\sigma}_{l\, \overline{l-1}} - q^{-1}
  \hat{\sigma}_{l\, \overline {l-1}} \hat{\sigma}_{bl},$&$
  \varepsilon_b > \varepsilon_l $\\ $\hat{\sigma}_{b \overline{l}} =
  \hat{\sigma}_{b\,l-1} \hat{\sigma}_{l-1 \, \overline{l}} - q^{-1}
  \hat{\sigma}_{l-1\,\overline{l}} \hat{\sigma}_{b\,l-1} ,$&$
  \varepsilon_b > \varepsilon_{l-1} $\\ $\hat{\sigma}_{la} =
  \hat{\sigma}_{l\, \overline{l-1}} \hat{\sigma}_{\overline {l-1}\, a}
  - q^{-1} \hat{\sigma}_{\overline{l-1}\, a} \hat{\sigma}_{l\,
  \overline{l-1}},$&$ \varepsilon_a < -\varepsilon_{l-1} $\\
  $\hat{\sigma}_{l-1\, a} = q^{-(\alpha_l, \varepsilon_a)}
  \hat{\sigma}_{l-1\, \overline{l}}\hat{\sigma}_{\overline{l}a}
  -q^{-1}\hat{\sigma}_{\overline{l}a}
  \hat{\sigma}_{l-1\,\overline{l}},$&$ \varepsilon_a <
  -\varepsilon_l$\\ $q^{(\alpha_l, \varepsilon_b)} \hat{\sigma}_{ba}
  \hat{\sigma}_{l-1\, \overline{l}}- q^{-(\alpha_l,\varepsilon_a)}
  \hat{\sigma}_{l-1\,\overline{l}} \hat{\sigma}_{ba} = 0,$&$
  \varepsilon_b > \varepsilon_a;\: a \neq l, l-1;\: b \neq
  \overline{l-1}, \overline{l} $ \\
\noalign{\smallskip}\hline
\end{tabular}
\end{table}
\begin{table}[h]
\caption{The relations for the operators $\hat{\sigma}_{ba}$ that hold
only for odd $m$} 
\label{odd} 
\centering
\begin{tabular}{|ll|} 
\hline\noalign{\smallskip}
\hspace{20mm} Relation & Conditions \\ 
\noalign{\smallskip}\hline\noalign{\smallskip}
$\hat{\sigma}_{b\,l+1} =
\hat{\sigma}_{bl} \hat{\sigma}_{l\,l+1} - q^{-1} \hat{\sigma}_{l\,l+1}
\hat{\sigma}_{bl},$&$ \varepsilon_b > \varepsilon_l$\\
$\hat{\sigma}_{b\, \overline{l}} = q^{(\alpha_l, \varepsilon_b)}
\hat{\sigma} _{b\,l+1} \hat{\sigma}_{l+1\,\overline{l}}
-\hat{\sigma}_{l+1\, \overline{l}} \hat{\sigma}_{b\, l+1},$&$
\varepsilon_b > 0$\\ $\hat{\sigma}_{la} = q^{-(\alpha_l,
\varepsilon_a)} \hat{\sigma}_{l\, l+1} \hat{\sigma}_{l+1\, a} -
\hat{\sigma}_{l+1\, a} \hat{\sigma}_{l\, l+1}, $&$ \varepsilon_a < 0
$\\ $\hat{\sigma}_{l+1\, a} = \hat{\sigma}_{l+1\, \overline{l}}
\hat{\sigma}_ {\overline{l}a} - q^{-1} \hat{\sigma}_{\overline{l}a}
\hat{\sigma}_{l+1\, \overline{l}}, $&$\varepsilon_a < -\varepsilon_l
$\\ $q^{(\alpha_l, \varepsilon_b)} \hat{\sigma}_{ba} \hat{\sigma}_{l\,
l+1} - q^{-(\alpha_l, \varepsilon_a)} \hat{\sigma}_{l\,l+1}
\hat{\sigma}_{ba} = 0,$ &$\varepsilon_b >\varepsilon_a;\: a \neq l,
l+1;\: b\neq l+1,\overline{l}$\\ 
\noalign{\smallskip}\hline
\end{tabular}
\end{table}

\newpage


\end{document}